\documentclass[a4paper,12pt]{article}
\usepackage{amssymb,latexsym,amsmath,amsthm,amscd,graphpap,dsfont}
\usepackage{amsmath}
\usepackage{amsfonts}
\usepackage{amssymb}
\usepackage[dvips]{epsfig}
\newtheorem{theorem}{Theorem}[section]
\newtheorem{prop}[theorem]{Proposition}

\newtheorem{lemma}[theorem]{Lemma}

\newtheorem{coro}[theorem]{Corollary}

\newenvironment{demo}{\noindent\textbf{Proof:}}{\hfill
$\square$\\ \medskip}
\newenvironment{rem}{\noindent\textbf{Remark:}}{\medskip}
\newenvironment{rems}{\noindent\textbf{Remarks:}}{\medskip}
\newcommand{\reff}[1]{(\ref{#1})}

\newcommand{\Cm}{\mathbb{C}}

\newcommand{\Nm}{\mathbb{N}}
\newcommand{\Rm}{\mathbb{R}}

\renewcommand{\AA}{\mathcal{A}}

\newcommand{\CC}{\mathcal{C}}

\newcommand{\EE}{\mathcal{E}}
\newcommand{\FF}{\mathcal{F}}

\newcommand{\LL}{\mathcal{L}}

\newcommand{\OO}{\mathcal{O}}

\newcommand{\dd}{\,\mathrm{d}}
\newcommand{\ul}{\mathrm{ul}}
\newcommand{\loc}{\mathrm{loc}}
\newcommand{\ess}{\mathrm{ess}}
\newcommand{\bu}{\mathrm{bu}}
\renewcommand{\epsilon}{\varepsilon}
\renewcommand{\Re}{\mathrm{Re}}

\numberwithin{equation}{section}

\newcounter{constantes}

\newdimen\texpscorrection
\newdimen\figcenter
\def\figurewithtex #1 #2 #3 #4 #5\cr{\null
  {\goodbreak\figcenter=\hsize\relax
  \advance\figcenter by -#4truecm
  \divide\figcenter by 2
  \begin{figure}[hbt]
  \vskip #3truecm\noindent\hskip\figcenter
  \includegraphics{#1}{\hskip\texpscorrection\input #2 }
  \vskip 0.8truecm{\baselineskip=0.8\baselineskip
  \noindent \vbox{\noindent {\footnotesize #5}}\par}
  \end{figure}}}
\def\point#1 #2 #3 {\rlap{\kern #1 truecm
\raise #2 truecm \hbox{#3}}}

\begin{document}

\title{\bf Global stability of travelling fronts for a damped wave 
equation with bistable nonlinearity} 

\author{Thierry GALLAY \& Romain JOLY\\
{Institut Fourier, UMR CNRS 5582}\\ 
{Universit\'e de Grenoble I}\\
{B.P. 74}\\
{38402 Saint-Martin-d'H\`eres, France}}

\date{October 3, 2007}

\maketitle

\begin{abstract}
We consider the damped wave equation $\alpha u_{tt}+u_t=u_{xx}-V'(u)$
on the whole real line, where $V$ is a bistable potential. This
equation has travelling front solutions of the form $u(x,t)=h(x-st)$
which describe a moving interface between two different steady states
of the system, one of which being the global minimum of $V$. We show
that, if the initial data are sufficiently close to the profile of a
front for large $|x|$, the solution of the damped wave equation
converges uniformly on $\Rm$ to a travelling front as $t \to +\infty$.
The proof of this global stability result is inspired by a recent work
of E.~Risler \cite{Ri} and relies on the fact that our system has
a Lyapunov function in any Galilean frame.\\[1mm]
{\bf Keywords:} travelling front, global stability, damped
wave equation, Lyapunov function.\\[1mm]
{\bf AMS classification codes (2000):} 35B35, 35B40, 37L15, 37L70.
\end{abstract}

\section{Introduction}
\label{intro}

The aim of this paper is to describe the long-time behavior of a large
class of solutions of the semilinear damped wave equation
\begin{equation}\label{eq-u}
  \alpha u_{tt}+ u_t\,=\, u_{xx} -  V'(u)~,
\end{equation}
where $\alpha > 0$ is a parameter, $V : \Rm \to \Rm$ is a smooth
bistable potential, and the unknown $u = u(x,t)$ is a real-valued
function of $x \in \Rm$ and $t \ge 0$. Equations of this form appear
in many different contexts, especially in physics and in biology. For
instance, Eq.~\eqref{eq-u} describes the continuum limit of an infinite
chain of coupled oscillators, the propagation of voltage along a
nonlinear transmission line \cite{DO}, and the evolution of an
interacting population if the spatial spread of the individuals is
modelled by a velocity jump process instead of the usual Brownian
motion \cite{Go,Ha3,Kac}.

As was already observed by several authors, the long-time asymptotics
of the solutions of the damped wave equation \eqref{eq-u} are quite
similar to those of the corresponding reaction-diffusion equation $u_t
= u_{xx} -V'(u)$. In particular, if $V'(u)$ vanishes rapidly enough as
$u \to 0$, the solutions of \eqref{eq-u} originating from small and
localized initial data converge as $t \to +\infty$ to the same
self-similar profiles as in the parabolic case
\cite{GR2,INZ,Kar,Ni1,Ni2}. The analogy persists for solutions with
nontrivial limits as $x \to \pm\infty$, in which case the long-time
asymptotics are often described by uniformly translating solutions of
the form $u(x,t)=h(x-st)$, which are usually called {\sl travelling
fronts}. Existence of such solutions for hyperbolic equations of the
form \eqref{eq-u} was first proved by Hadeler \cite{Ha1,Ha2}, and a
few stability results were subsequently obtained by Gallay \& Raugel
\cite{Ga,GR1,GR3,GR4}.

While local stability is an important theoretical issue, in the
applications one is often interested in {\em global convergence
  results} which ensure that, for a large class of initial data with a
prescribed behavior at infinity, the solutions approach travelling
fronts as $t \to +\infty$. For the scalar parabolic equation $u_t =
u_{xx} - V'(u)$, such results were obtained by Kolmogorov, Petrovski
\& Piskunov \cite{KPP}, by Kanel \cite{kanel1,kanel2}, and by Fife \&
McLeod \cite{FM1,FM2} under various assumptions on the potential.  All
the proofs use in an essential way comparison theorems based on the
maximum principle. These techniques are very powerful to obtain global
information on the solutions, and were also successfully applied to
monotone parabolic systems \cite{VVV,RTV} and to parabolic equations
on infinite cylinders \cite{roque1,roque2}.

However, unlike its parabolic counterpart, the damped wave equation
\eqref{eq-u} has no maximum principle in general. More precisely,
solutions of \eqref{eq-u} taking their values in some interval $I
\subset \Rm$ obey a comparison principle only if
\begin{equation}\label{pm}
   4 \alpha \sup_{u \in I} V''(u) \,\le\, 1~,
\end{equation}
see \cite{PW} or \cite[Appendix~A]{GR1}. In physical terms, this
condition means that the relaxation time $\alpha$ is small compared to
the period of the nonlinear oscillations. In particular, if $I$ is a
neighborhood of a local minimum $\bar u$ of $V$, inequality \eqref{pm}
implies that the linear oscillator $\alpha u_{tt} + u_t + V''(\bar u)u
= 0$ is strongly damped, so that no oscillations occur. It was shown
in \cite{GR1,GR3} that the travelling fronts of \eqref{eq-u} with a
monostable nonlinearity are stable against large perturbations
provided that the parameter $\alpha$ is sufficiently small so that the
strong damping condition \eqref{pm} holds for the solutions under
consideration.  In other words, the basin of attraction of the
hyperbolic travelling fronts becomes arbitrarily large as $\alpha \to
0$, but if $\alpha$ is not assumed to be small there is no hope to use
``parabolic'' methods to obtain global stability results for the
travelling fronts of the damped wave equation \eqref{eq-u}.

Recently, however, a different approach to the stability of travelling
fronts has been developped by Risler \cite{GR5,Ri}. The new method is
purely variational and is therefore restricted to systems that possess
a gradient structure, but its main interest lies in the fact that it
does not rely on the maximum principle. The power of this approach is
demonstrated in the pioneering work \cite{Ri} where global convergence
results are obtained for the non-monotone reaction-diffusion system
$u_t = u_{xx} - \nabla V(u)$, with $u \in \Rm^n$ and $V : \Rm^n \to
\Rm$. The aim of the present article is to show that Risler's method
can be adapted to the damped hyperbolic equation \eqref{eq-u} and
allows in this context to prove global convergence results {\em without
any smallness assumption} on the parameter $\alpha$. This will also
be an opportunity to present the main arguments of \cite{Ri} in an
alternative way, although some important ingredients of our proof
have no counterpart in Risler's work. 

\medskip
Before stating our theorem, we need to specify the assumptions we make
on the nonlinearity in \eqref{eq-u}. We suppose that $V \in
\CC^3(\Rm)$, and that there exist positive constants $a$, $b$ such that
\begin{equation}\label{coercive}
   u V'(u) \,\ge\, a u^2 - b~, \quad \hbox{for all } u \in \Rm~.
\end{equation}
In particular, $V(u) \to +\infty$ as $|u| \to \infty$. We also 
assume
\begin{eqnarray}\label{u=0}
   &&V(0) \,=\, 0~, \quad V'(0) \,=\, 0~, \quad V''(0) \,>\, 0~, \\
   \label{u=1}
   &&V(1) \,<\, 0~, \quad V'(1) \,=\, 0~, \quad V''(1) \,>\, 0~. 
\end{eqnarray}
Finally we suppose that, except for $V(0)$ and $V(1)$, all 
critical values of $V$ are positive: 
\begin{equation}\label{critval}
  \Bigl\{u \in \Rm \,\Big|\, V'(u) = 0\,,~V(u) \le 0\Bigr\} 
  \,=\, \{0\,;\,1\}~.
\end{equation}
In other words $V$ is a smooth, strictly coercive function 
which reaches its global minimum at $u = 1$ and has in addition 
a local minimum at $u = 0$. We call $V$ a {\em bistable} 
potential because both $u = 0$ and $u = 1$ are stable equilibria
of the one-dimensional dynamical system $\dot u = -V'(u)$. 
The simplest example of such a potential is represented in 
Fig.~1. Note however that $V$ is allowed to have positive critical 
values, including local minima. 

\figurewithtex 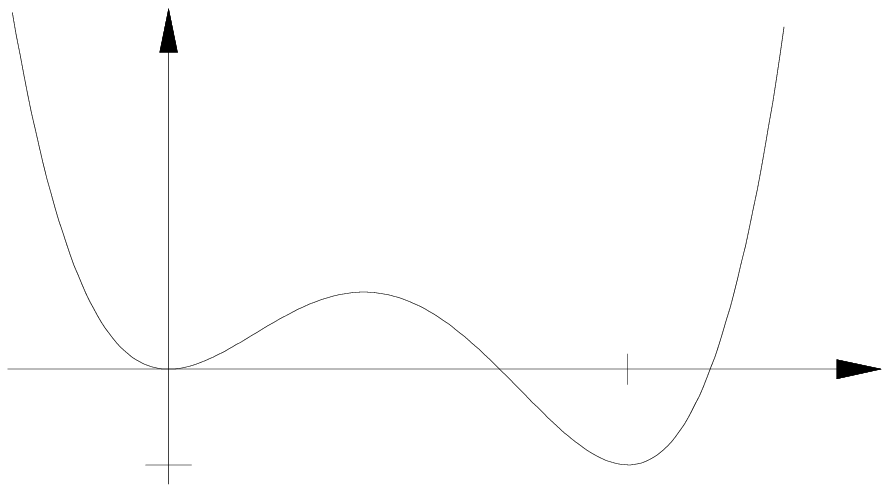 fig1.tex 5.000 9.000
{\bf Fig.~1:} The simplest example of a potential $V$ satisfying 
assumptions \eqref{coercive}--\eqref{critval}.\cr

Under assumptions \eqref{u=0}--\eqref{critval}, it is well-known that
the parabolic equation $u_t = u_{xx} -V'(u)$ has a family of 
travelling fronts of the form $u(x,t) = h(x-c_*t-x_0)$ connecting 
the stable equilibria $u = 1$ and $u = 0$, see e.g. \cite{AW}. 
More precisely, there exists a unique speed $c_* > 0$ such that 
the boundary value problem
\begin{equation}\label{hdef}
\left\{\begin{array}{l}
  h''(y) + c_* h'(y) - V'(h(y)) \,=\, 0~, \quad y \in \Rm~,\\
  h(-\infty) = 1~, \quad h(+\infty) = 0~,
\end{array}\right.
\end{equation}
has a solution $h : \Rm \to (0,1)$, in which case the profile 
$h$ itself is unique up to a translation. Moreover $h \in
\CC^4(\Rm)$, $h'(y) < 0$ for all $y \in \Rm$, and $h(y)$ 
converges exponentially toward its limits as $y \to \pm\infty$. 
As was observed in \cite{GR1,Ha1}, for any $\alpha > 0$ the
damped hyperbolic equation \eqref{eq-u} has a corresponding 
family of travelling fronts given by
\begin{equation}\label{TWdef}
  u(x,t) \,=\, h(\sqrt{1+\alpha c_*^2}\,x - c_*t - x_0)~, 
  \quad x_0 \in \Rm~.
\end{equation}
Remark that the actual speed of these waves is not $c_*$, 
but $s_* = c_*/\sqrt{1+\alpha c_*^2}$. In particular $s_*$ 
is smaller than $1/\sqrt{\alpha}$ (the slope of the characteristics 
of Eq.~\eqref{eq-u}), which means that the travelling fronts 
\eqref{TWdef} are always ``subsonic''. In what follows we shall
refer to $c_*$ as the ``parabolic speed'' to distinguish it 
from the physical speed $s_*$. 

We are now in position to state our main result:

\begin{theorem}\label{mainthm}
Let $\alpha > 0$ and let $V \in \CC^3(\Rm)$ satisfy 
\eqref{coercive}--\eqref{critval} above. Then there exist positive 
constants $\delta$ and $\nu$ such that, for all initial data 
$(u_0,u_1) \in H^1_\ul(\Rm) \times L^2_\ul(\Rm)$ satisfying
\begin{eqnarray}\label{cond-}
  \limsup_{\xi \to -\infty} \int_\xi^{\xi+1} \Bigl((u_0(x)-1)^2 + {u_0'(x)}^2
  + u_1(x)^2\Bigr) \dd x &\le& \delta~,\\ \label{cond+}
  \limsup_{\xi \to +\infty} \int_\xi^{\xi+1} \Bigl(u_0(x)^2 + {u_0'(x)}^2
  + u_1(x)^2\Bigr) \dd x &\le& \delta~,
\end{eqnarray}
equation \eqref{eq-u} has a unique global solution (for positive 
times) such that $u(\cdot,0) = u_0$, $u_t(\cdot,0) = u_1$. 
Moreover, there exists $x_0 \in \Rm$ such that
\begin{equation}\label{convres}
  \sup_{x \in \Rm}\Big| u(x,t) - h(\sqrt{1+\alpha c_*^2}\,x 
  - c_*t - x_0)\Big| \,=\, \OO(e^{-\nu t})~, \quad 
  \hbox{as } t \to +\infty~. 
\end{equation}
\end{theorem}

\medskip\noindent{\bf Remarks:}\\
{\bf 1.} Loosely speaking Theorem~\ref{mainthm} says that, if the
initial data $(u_0,u_1)$ are close enough to the global equilibrium
$(1,0)$ as $x \to -\infty$ and to the local equilibrium $(0,0)$ as $x
\to +\infty$, the solution $u(x,t)$ of \eqref{eq-u} converges
uniformly in space and exponentially fast in time toward a member of
the family of travelling fronts \eqref{TWdef}.  In particular, any
solution which looks roughly like a travelling front at initial time
will eventually approach a suitable translate of that front. It should
be noted, however, that our result does not give any constructive
estimate of the time needed to reach the asymptotic regime described
by \eqref{convres}.  Depending on the shape of the potential and of
the initial data, very long transients can occur before the solution
actually converges to a travelling front.\\[1mm]
{\bf 2.} The definition of the uniformly local Lebesgue space
$L^2_\ul(\Rm)$ and the uniformly local Sobolev space $H^1_\ul(\Rm)$
will be recalled at the beginning of Section~\ref{sec2}. These spaces
provide a very convenient framework to study infinite-energy solutions
of the hyperbolic equation \eqref{eq-u}, but their knowledge is not
necessary to understand the meaning of Theorem~\ref{mainthm}. In a
first reading one can assume, for instance, that $u_0'$ and $u_1$ are
bounded and uniformly continuous functions, in which case assumptions
\eqref{cond-}, \eqref{cond+} can be replaced by 
\[ 
  \limsup_{x \to -\infty} (|u_0(x)-1| + |u_0'(x)| + |u_1(x)|) 
  \,\le\, \delta~, \quad \limsup_{x \to +\infty} (|u_0(x)| + 
  |u_0'(x)| + |u_1(x)|) \,\le\, \delta~.
\]
Also, to simplify the presentation, we have expressed our convergence
result \eqref{convres} in the uniform norm, but the proof will 
show that the solution $u(x,t)$ of \eqref{eq-u} converges to a 
travelling front in the uniformly local energy space $H^1_\ul(\Rm) 
\times L^2_\ul(\Rm)$, see \eqref{convul} below. \\[1mm]
{\bf 3.} The convergence rate $\nu$ in \eqref{convres} is related to
the spectral gap of the linearization of \eqref{eq-u} at the
travelling front. As is shown in Section~\ref{sec9}, we can take 
$\nu = \OO(1)$ in the parabolic limit $\alpha \to 0$, whereas 
$\nu = \OO(1/\alpha)$ as $\alpha \to +\infty$. On the other hand, 
the parameter $\delta$ in \eqref{cond-}, \eqref{cond+} must be chosen 
small enough so that the following two conditions are satisfied. 
First, the initial data $(u_0,u_1)$ should lie in the local basin 
of attraction of the steady state $(0,0)$ for large positive $x$, 
and in the basin of $(1,0)$ for large negative $x$. Next, the 
energy integral
\[
  \int_{-\infty}^0 e^{cx} \Bigl(\frac{\alpha}{2}u_1(x)^2 + \frac12
  u_0'(x)^2 + V(u_0(x))\Bigr)\dd x~,
\]
which is well-defined for any $c > 0$, should diverge to $-\infty$
as $c \to 0$. The second condition is an essential ingredient of
our variational proof, but we do not know if the conclusion of 
Theorem~\ref{mainthm} still holds without such an assumption. 

\medskip
The proof of Theorem~\ref{mainthm} relies on the fact that
Eq.~\eqref{eq-u} has, at least formally, a whole family of 
{\em Lyapunov functions}. To see this, let $u(x,t)$ be a solution 
of \eqref{eq-u} whose initial data satisfy \eqref{cond-}, \eqref{cond+}. 
Given any $c \ge 0$ we go to a uniformly translating frame by 
setting
\begin{equation}\label{def-uc}
  u(x,t) \,=\, u_c(\sqrt{1+\alpha c^2}\,x - ct\,,\,t)~, 
  \quad \hbox{or} \quad 
  u_c(y,t) \,=\, u\Bigl(\frac{x+ct}{\sqrt{1+\alpha c^2}}\,,\,t
  \Bigr)~. 
\end{equation}
The new function $u_c(y,t)$ is then a solution of the modified 
equation
\begin{equation}\label{eq-uc}
  \alpha \ddot u_c + \dot u_c - 2\alpha c \dot u'_c \,=\, 
  u''_c + c u'_c - V'(u_c)~,
\end{equation}
where $\dot u_c(y,t) \equiv \partial_t u_c(y,t)$ and  
$u'_c(y,t) \equiv \partial_y u_c(y,t)$. If we now introduce 
the energy function
\begin{equation}\label{Ecdef}
  E_c(t) \,=\, \int_\Rm e^{cy} \Bigl(\frac{\alpha}{2}
  (\dot u_c(y,t))^2 + \frac12(u'_c(y,t))^2 
  + V(u_c(y,t))\Bigr)\dd y~,
\end{equation}
a direct calculation shows that 
\begin{equation}\label{dotEc}
  E_c'(t) \,=\, -(1+\alpha c^2)\int_\Rm e^{cy} (\dot u_c(y,t))^2 
  \dd y \,\le\, 0~.
\end{equation}
In other words, Eq.~\eqref{eq-u} possess (at least formally) a
continuous family of non-equivalent Lyapunov functions, indexed by the
parabolic speed $c \ge 0$. In the parabolic case $\alpha = 0$, it is
shown in \cite{GR5} that this rich Lyapunov structure is sufficient to
prove the convergence \eqref{convres} if we restrict ourselves to
solutions which decay sufficiently rapidly to zero as $x \to +\infty$, 
and we believe that the approach of \cite{GR5} works in the hyperbolic 
case too. 

However, it is important to realize that the solutions we consider in
Theorem~\ref{mainthm} are only supposed to be small for large positive
$x$, and do not necessarily converge to zero as $x \to +\infty$.
Under these assumptions the integral in \eqref{Ecdef}, which contains
the exponentially growing factor $e^{cy}$, is usually divergent at
$+\infty$, so that the Lyapunov function $E_c$ is certainly not
well-defined. This is a technical problem which seriously complicates
the analysis. To overcome this difficulty, a possibility is to
truncate the exponential factor $e^{cy}$ in \eqref{Ecdef} to make it
integrable over $\Rm$, see \cite{FM1}, \cite{Ri}. We choose here
another solution which consists in decomposing the solution $u(x,t)$
into a principal part $v(x,t)$ which is compactly supported to the
right, and a small remainder $r(x,t)$ which decays exponentially to
zero as $t \to +\infty$. The idea is then to study the approximate
Lyapunov function defined by \eqref{Ecdef} with $u_c(x,t)$ replaced by
$v_c(x,t)$, see Section~\ref{sec4} for more details.

As was already mentioned, the proof of Theorem~\ref{mainthm} closely
follows the previous work \cite{Ri} which deals with gradient
reaction-diffusion systems of the form $u_t = u_{xx} - \nabla V(u)$.
There are, however, significant differences that we want to emphasize.
First, the evolution defined by the damped hyperbolic equation
\eqref{eq-u} is not regularizing in finite time, but only
asymptotically as $t \to +\infty$. As a consequence, the compactness
arguments which play an essential role in the proof become slightly
more delicate in the hyperbolic case. On the other hand, the solutions
of \eqref{eq-u} have a {\em finite speed of propagation}, a property
which has no parabolic analog. Although this is not an essential
ingredient of the proof, we shall take advantage of this fact here and
there to get a priori estimates on the solutions of \eqref{eq-u}.
Finally, an important property of the {\em scalar} equation
\eqref{eq-u} is that the associated elliptic problem \eqref{hdef} has
a unique solution $(h,c_*)$, and that the corresponding travelling
front is a stable solution of \eqref{eq-u}. This is no longer true
for the systems considered in \cite{Ri}, in which several stable or
unstable fronts with different speeds may connect the same pair of
equilibria. In this more general situation, without additional
assumptions one can only show that the solution $u(x,t)$ approaches
as $t \to \infty$ the family of all travelling fronts with a 
given speed. 

Besides these natural differences due to the properties of 
Eq.~\eqref{eq-u}, we also made technical choices in our proof 
which substantially differ from \cite{Ri}. As was already mentioned, 
the most important one is that we give a meaning to the Lyapunov
function $E_c$ by decomposing the solution $u(x,t)$, and not 
by truncating the exponential weight $e^{cy}$. The main avantage
of this approach is that the behavior of the energy is then easier 
to control. However, new arguments are required which have 
no counterpart in \cite{Ri} or \cite{GR5}. This is the case in 
particular of Section~\ref{sec6}, which is the main technical 
step in our proof. 

\medskip The rest of this paper is organized as follows. In
Section~\ref{sec2}, we briefly present the uniformly local spaces and
we study the Cauchy problem for Eq.~\eqref{eq-u} in this framework. In
Section~\ref{sec3}, we prove the persistence of the boundary
conditions \eqref{cond-}, \eqref{cond+} and we decompose the solution
of \eqref{eq-u} as $u(x,t) = v(x,t) + r(x,t)$, where $v$ is compactly
supported to the right and $r$ decays exponentially as $t \to
+\infty$. We also introduce the invasion point $\overline x(t)$ which
tracks the position of the moving interface. The core of the proof
starts in Section~\ref{sec4}, where we control the behavior of the
energy $E_c$ in a frame moving at constant speed $s = c/\sqrt{1+
\alpha c^2}$. These estimates are used in Section~\ref{sec5} 
to prove that the average speed $\overline x(t)/t$ converges to 
a limit $s_\infty\in(0,1/\sqrt{\alpha})$ as $t \to +\infty$. 
The main technical step is Section~\ref{sec6}, where we show that 
the energy stays uniformly bounded in a frame following the invasion 
point, see Proposition~\ref{prop-control-energy} for a precise statement. 
This allows us to prove in Section~\ref{sec7} that the solution 
$u(x,t)$ converges as $t \to +\infty$ to a travelling front 
uniformly in any interval of the form $(\overline x(t)-L,+\infty)$.
The proof of Theorem~\ref{mainthm} is then completed in two steps. 
In Section~\ref{sec8}, we use an energy estimate in the laboratory 
frame to show that the solution $u(x,t)$ converges uniformly 
on $\Rm$ to a travelling front, at least for a sequence of times.
Finally, the local stability result established in Section~\ref{sec9}
gives the convergence for all times and the exponential rate 
in \eqref{convres}. 

\medskip\noindent{\bf Notations.} The symbols $K_0, K_1, \dots$
denote our main constants, which will be used throughout the paper. 
In contrast, the local constants $C_0, C_1, \dots$ will change 
from a section to another. We also denote by $C$ a positive constant
which may change from place to place.

\medskip\noindent{\bf Acknowledgements.} As is emphasized in the text
our approach is essentially based on ideas and techniques introduced
by Emmanuel Risler, to whom we are also indebted for many fruitful 
discussions. The work of Th.G was partially supported by the French 
Ministry of Research through grant ACI JC 1039.


\section{Global existence and asymptotic compactness}
\label{sec2}

In this section, we prove that the Cauchy problem for Eq.~\eqref{eq-u}
is globally well-posed for positive times in the uniformly local
energy space $X = H^1_\ul(\Rm) \times L^2_\ul(\Rm)$. We first recall 
the definitions of the uniformly local Sobolev spaces which provide
a natural framework for the study of partial differential equations
on unbounded domains, see e.g. \cite{ARCD,EZ,Fe,GS,GV,Ka,MT,MS1,MS2}. 

For any $u \in L^2_\loc(\Rm)$ we denote
\begin{equation}\label{L2ulnorm}
   \|u\|_{L^2_\ul} \,=\, \sup_{\xi \in \Rm} \Bigl(\int_\xi^{\xi+1} 
   |u(x)|^2 \dd x \Bigr)^{1/2} \,=\, \sup_{\xi \in \Rm} 
   \|u\|_{L^2([\xi,\xi+1])} \,\le\, \infty~.
\end{equation}
The uniformly local Lebesgue space is defined as
\begin{equation}\label{L2uldef}
   L^2_\ul(\Rm) \,=\, \Bigl\{u \in L^2_\loc(\Rm) \,\Big|\, 
   \|u\|_{L^2_\ul} < \infty\,,~\lim_{\xi \to 0}
   \|T_\xi u - u\|_{L^2_\ul} = 0\Bigr\}~,
\end{equation}
where $T_\xi$ denotes the translation operator: $(T_\xi u)(x) = 
u(x-\xi)$. In a similar way, for any $k \in \Nm$, we introduce 
the uniformly local Sobolev space
\begin{equation}\label{Hkuldef}
  H^k_\ul(\Rm) \,=\, \Bigl\{u \in H^k_\loc(\Rm) \,\Big|\,
  \partial^j u \in L^2_\ul(\Rm) \hbox{ for } j = 0,1,2,\dots,k
  \Bigr\}~,
\end{equation}
which is equipped with the natural norm
\[
  \|u\|_{H^k_\ul} \,=\, \Bigl(\sum_{j=0}^k \|\partial^j 
  u\|_{L^2_\ul}^2 \Bigr)^{1/2}~.
\]
It is easy to verify that $H^k_\ul(\Rm)$ is a Banach space, which is 
however neither reflexive nor separable. If $\CC^k_\bu(\Rm)$ denotes
the Banach space of all $u \in \CC^k(\Rm)$ such that $\partial^j u$
is bounded and uniformly continuous for $j = 0,\dots,k$, we have 
the continuous inclusions 
\[
  \CC^k_\bu(\Rm) \,\hookrightarrow\, H^k_\ul(\Rm) 
  \,\hookrightarrow\, \CC^{k-1}_\bu(\Rm)~.
\]
In particular $H^1_\ul(\Rm) \hookrightarrow \CC^0_\bu(\Rm)$ and 
$\|u\|_{L^\infty}^2 \le 2\|u\|_{H^1_\ul}^2$ for all $u \in H^1_\ul$. 
Note also that $H^k_\ul(\Rm)$ is an algebra for any $k \ge 1$, i.e. 
$\|uv\|_{H^k_\ul} \le C \|u\|_{H^k_\ul} \|v\|_{H^k_\ul}$ for all 
$u,v \in H^k_\ul(\Rm)$. Finally the space $\CC^\infty_\bu(\Rm)$ is 
dense in $H^k_\ul(\Rm)$ for any $k \in \Nm$. 

\smallskip
\begin{rem} Some authors do not include in the definition of the
uniformly local $L^2$ space the assumption that $\xi \mapsto T_\xi u$ 
is continuous for any $u \in L^2_\ul(\Rm)$. The resulting uniformly 
local Sobolev spaces are of course larger, but also less convenient
from a functional-analytic point of view. In particular, one looses
the property that $H^{k+1}_\ul(\Rm)$ is dense in $H^k_\ul(\Rm)$. As 
we shall see, the definitions \eqref{L2uldef}, \eqref{Hkuldef} 
guarantee that the damped wave equation \eqref{eq-u} defines 
a continuous evolution in $H^1_\ul(\Rm) \times L^2_\ul(\Rm)$.
\end{rem}

Let $X = H^1_\ul(\Rm) \times L^2_\ul(\Rm)$ and $Y = H^2_\ul(\Rm)
\times H^1_\ul(\Rm)$. The main result of this section is:

\begin{prop}\label{propexist} For all initial data $(u_0,u_1) \in X$, 
Eq.~\eqref{eq-u} has a unique global solution $u \in \CC^0([0,+\infty),
H^1_\ul(\Rm)) \cap \CC^1([0,+\infty),L^2_\ul(\Rm))$ satisfying 
$u(0) = u_0$, $u_t(0) = u_1$. This solution depends continuously 
on the initial data, uniformly in time on compact intervals. 
Moreover, there exists $K_\infty > 0$ (depending only on $\alpha$ and 
$V$) such that
\begin{equation}\label{attball}
  \limsup_{t \to +\infty}\,(\|u(\cdot,t)\|_{H^1_\ul}^2 +
  \|u_t(\cdot,t)\|_{L^2_\ul}^2) \,\le\, K_\infty~.
\end{equation}
\end{prop}

\begin{demo}
Setting $w = (u,u_t)$, we rewrite \eqref{eq-u} as a first order 
evolution equation 
\begin{equation}\label{eq-w}
  w_t \,=\, A w + F(w)~,
\end{equation}
where
\begin{equation}\label{AFdef}
  A \,=\, \frac{1}\alpha \begin{pmatrix} 0 & \alpha \\ 
  \partial_x^2-1 & -1 \end{pmatrix}~, \quad \hbox{and} \quad 
  F(w) \,=\, \frac{1}{\alpha} \begin{pmatrix} 0 \\ 
  u - V'(u) \end{pmatrix}~.
\end{equation}
Using d'Alembert's formula for the solution of the wave equation 
$\alpha u_{tt} = u_{xx}$, it is straightforward to verify that the 
linear operator $A_0$ on $X$ defined by 
\[
  D(A_0) \,=\, Y~, \quad A_0 \,=\, \frac{1}\alpha \begin{pmatrix} 
  0 & \alpha \\ \partial_x^2 & 0 \end{pmatrix}~,
\]
is the generator of a strongly continuous group of bounded linear
operators in $X$. The same is true for the linear operator $A$, which
is a bounded perturbation of $A_0$, see \cite[Section~3.1]{Pa}.  On
the other hand, as $V \in \CC^3(\Rm)$ and $H^1_\ul(\Rm)
\hookrightarrow \CC^0_\bu(\Rm)$, it is clear that the nonlinearity $F$
maps $X$ into $Y$, and that $F$ is Lipschitz continuous on any bounded
set $B \subset X$.  Thus a classical argument shows that the Cauchy
problem for \eqref{eq-w} is locally well-posed in $X$, see
\cite[Section~6.1]{Pa} or \cite[Section~7.2]{GS}. More precisely, for
any $r > 0$, there exists $T(r) > 0$ such that, for all initial data
$w_0 \in X$ with $\|w_0\|_X \le r$, Eq.~\eqref{eq-w} has a unique (mild)
solution $w \in \CC^0([0,T],X)$ satisfying $w(0) = w_0$. This solution
$w(t)$ depends continuously on the initial data $w_0$ in $X$,
uniformly for all $t \in [0,T]$. Moreover, if $w_0 \in Y$, then 
$w \in \CC^1([0,T],X) \cap \CC^0([0,T],Y)$ is a classical solution of
\eqref{eq-w}. To prove Proposition~\ref{propexist}, it remains to show
that all solutions of \eqref{eq-w} stay bounded for positive times
(hence can be extended to global solutions), and are eventually
contained in an attracting ball whose radius is independent of the
initial data.

Assume that $w = (u,u_t) \in \CC^0([0,T],X)$ is a solution of 
\eqref{eq-w}. Let $\rho(x) = \exp(-\kappa|x|)$, where $\kappa > 0$ 
is small enough so that $2\sqrt{\alpha}\kappa \le 1$ and 
$\kappa^2 \le a$, with $a > 0$ as in \eqref{coercive}. For any 
$\xi \in \Rm$ and all $t \in [0,T]$, we define
\begin{eqnarray}\label{EEdef}
  \EE(\xi,t) &=& \int_\Rm (T_\xi \rho)(x)\Bigl(\alpha^2 u_t^2 
    + \alpha u_x^2 + 2\alpha \overline{V}(u) + \frac12 u^2 
    + \alpha u u_t\Bigr)(x,t)\dd x~, \\ \label{FFdef}
  \FF(\xi,t) &=& \int_\Rm (T_\xi \rho)(x)\Bigl(\alpha u_t^2 + u_x^2 
    + a u^2\Bigr)(x,t)\dd x~,
\end{eqnarray}
where $\overline{V}(u) = V(u)-V(1) \ge 0$ and $(T_\xi\rho)(x) = 
\rho(x-\xi)$. We also denote
\begin{equation}\label{Mdef}
  M(t) \,=\, \sup_{\xi \in \Rm} \EE(\xi,t)~, \quad t \in [0,T]~.
\end{equation}
Since $u(\cdot,t) \in H^1_\ul(\Rm)$ and $u_t(\cdot,t) \in
L^2_\ul(\Rm)$, it is clear that $M(t) < \infty$ for all $t \in [0,T]$.
Moreover, as $\overline{V}(u) \ge 0$ and $|\alpha u u_t| \le
\frac{3\alpha^2}4 u_t^2 + \frac13 u^2$, we have
\[
  \EE(\xi,t) \,\ge\, \int_\Rm (T_\xi \rho)(x)\Bigl(\frac{\alpha^2}4 
  u_t^2 + \alpha u_x^2 + \frac16 u^2\Bigr)(x,t)\dd x~.
\]
Taking in both sides the supremum over $\xi \in \Rm$ and using the
definitions \eqref{L2ulnorm}--\eqref{Hkuldef}, we see that there
exists $C_1 > 0$ (depending only on $\alpha$) such that 
$\|w(\cdot,t)\|_X^2 \equiv \|u(\cdot,t)\|_{H^1_\ul}^2 + 
\|u_t(\cdot,t)\|_{L^2_\ul}^2 \le C_1 M(t)$. 

On the other hand, differentiating $\EE(\xi,t)$ with respect to 
time, we find
\[
  \partial_t \EE(\xi,t) \,=\, -\int_\Rm 
    (T_\xi \rho)(\alpha u_t^2 + u_x^2 + uV'(u))\dd x 
  -\int_\Rm (T_\xi \rho)'(uu_x + 2\alpha u_x u_t)\dd x~.
\]
To estimate the last integral, we observe that 
\[
  -\int_\Rm (T_\xi \rho)'uu_x\dd x \,=\, \frac{\kappa^2}2  
  \int_\Rm (T_\xi \rho) u^2 \dd x - \kappa u(\xi)^2 \,\le\, 
  \frac{\kappa^2}2 \int_\Rm (T_\xi \rho)u^2\dd x~,
\]
and
\[
  \Bigl| \int_\Rm (T_\xi \rho)' 2\alpha u_x u_t\dd x \Bigr|
  \,\le\, 2\alpha\kappa \int_\Rm (T_\xi \rho)|u_x u_t| \dd x 
  \,\le\, \sqrt{\alpha}\kappa \int_\Rm (T_\xi \rho)
  (\alpha u_t^2 + u_x^2)\dd x~. \\
\]
Using \eqref{coercive} together with our assumptions on $\kappa$, 
we arrive at
\begin{equation}\label{dotE0}
  \partial_t \EE(\xi,t) \,\le\,  -\frac12 
  \int_\Rm (T_\xi \rho)(\alpha u_t^2 + u_x^2 + au^2)\dd x 
  + \frac{2b}{\kappa} \,=\, -\frac12 \FF(\xi,t) + \frac{2b}{\kappa}~. 
\end{equation}
This differential inequality implies that the quantity $M(t)$ defined
in \eqref{Mdef} is a decreasing function of time as long as it stays 
above a certain threshold. More precisely, we have:

\begin{lemma}\label{attlem}
There exists $C_2 > 0$ (depending only on $\alpha, V$) such that, 
if $M(t) \ge C_2$ for some $t \in [0,T]$ and $\EE(\xi,t)\ge M(t)-1$ 
for some $\xi \in \Rm$, then $\partial_t \EE(\xi,t) \le -1$.
\end{lemma}

Assuming this result to be true, we now conclude the proof of
Proposition~\ref{propexist}. It follows readily from
Lemma~\ref{attlem} that $M(t) \le \max(C_2,M(0)-t)$ for all $t \in
[0,T]$, an estimate which holds for any solution $w \in
\CC^0([0,T],X)$ of \eqref{eq-w}. This shows that any solution of
\eqref{eq-w} stays bounded in $X$ for positive times (hence can be
extended to a global solution), and that \eqref{attball} holds with
$K_\infty = C_1 C_2$.
\end{demo}

\noindent{\bf Proof of Lemma~\ref{attlem}.}
We use the same notations as in the proof of
Proposition~\ref{propexist}. Let $C_3 = 2(1+2b/\kappa)$, and take $L >
0$ large enough so that $e^{\kappa L} \ge 3$. Fix also $\xi \in \Rm$
and $t \in [0,T]$. If $\FF(\xi,t) \ge C_3$, then $\partial_t
\EE(\xi,t) \le -1$ by \eqref{dotE0}.  On the other hand, if
$\FF(\xi,t) \le C_3$, there exists $C_4 > 0$ (depending on $\alpha$,
$V$, $L$, and $C_3$) such that
\begin{equation}\label{C4bdd}
  \int_{\xi-L}^{\xi+L} (T_\xi \rho)(x) \,e(u,u_x,u_t)(x,t) \dd x 
  \,\le\, C_4~,
\end{equation}
where $e(u,u_x,u_t) = \alpha^2 u_t^2 + \alpha u_x^2 + 2\alpha 
\overline{V}(u) + \frac12 u^2 + \alpha u u_t \ge 0$. Inequality 
\eqref{C4bdd} holds because $\FF(\xi,t)$ controls the norm of 
$(u,u_t)$ in $H^1([\xi-L,\xi+L]) \times L^2([\xi-L,\xi+L])$. 
As a consequence of \eqref{EEdef}, \eqref{C4bdd} at least one 
of the following inequalities holds:
\begin{eqnarray}\label{altern}
  \int_{\xi+L}^\infty (T_\xi \rho)(x) \,e(u,u_x,u_t)(x,t) \dd x &\ge& 
    \frac12(\EE(\xi,t)-C_4)~, \quad\hbox{or}\quad \\ \nonumber 
  \int_{-\infty}^{\xi-L} (T_\xi \rho)(x) \,e(u,u_x,u_t)(x,t) \dd x &\ge& 
    \frac12(\EE(\xi,t)-C_4)~.
\end{eqnarray}
Suppose for instance that the first inequality in \eqref{altern} 
holds. Then 
\begin{eqnarray*}
  \EE(\xi+L,t) &\ge& \int_{\xi+L}^\infty (T_{\xi+L}\rho)(x) 
    \,e(u,u_x,u_t)(x,t) \dd x \\ 
  &\ge& 3 \int_{\xi+L}^\infty (T_\xi \rho)(x) \,e(u,u_x,u_t)(x,t) 
    \dd x \,\ge\, \frac32(\EE(\xi,t)-C_4)~,
\end{eqnarray*}
because $(T_{\xi+L}\rho)(x) \ge 3 (T_\xi \rho)(x)$ for all $x \ge 
\xi+L$, by assumption on $L$. Using a similar argument in the other 
case we conclude that, if $\FF(\xi,t) \le C_3$, then
\begin{equation}\label{notmax}
  \max\Bigl(\EE(\xi+L,t), \EE(\xi-L,t)\Bigr) \,\ge\, 
  \frac32(\EE(\xi,t)-C_4)~.
\end{equation}

Now, fix $C_5 > 3(C_4+1)$. If $M(t) \ge C_5$ and $\EE(\xi,t) \ge
M(t)-1$, we claim that $\FF(\xi,t) > C_3$, so that $\partial_t 
\EE(\xi,t) \le -1$ by \eqref{dotE0}. Indeed, if $\FF(\xi,t) \le 
C_3$, it follows from \eqref{notmax} that
\[
  M(t) \,\ge\, \max\Bigl(\EE(\xi+L,t), \EE(\xi-L,t)\Bigr) \,\ge\, 
  \frac32(M(t)-1-C_4)~, 
\]
which contradicts the assumption that $M(t) \ge C_5$. \hfill $\square$

\medskip
\begin{rem}
The proof of Proposition~\ref{propexist} can be simplified if we 
assume, in addition to \eqref{coercive}, that $uV'(u) \ge a'V(u)
-b'$ for some positive constants $a', b'$, but Lemma~\ref{attlem}
allows us to avoid this unnecessary assumption.
\end{rem}

\medskip
To conclude this section, we show that the solutions of \eqref{eq-u}
given by Proposition~\ref{propexist} are locally asymptotically 
compact, in the following sense: 

\begin{prop}\label{propcompact} Let $u \in \CC^0([0,+\infty),
H^1_\ul(\Rm)) \cap \CC^1([0,+\infty),L^2_\ul(\Rm))$ be a solution
of \eqref{eq-u}, and let $\{(x_n,t_n)\}_{n\in\Nm}$ be a sequence
in $\Rm \times \Rm_+$ such that $t_n \to +\infty$ as $n \to \infty$.
Then there exists a subsequence, still denoted $(x_n,t_n)$,  
and a solution $\bar u \in \CC^0(\Rm,H^1_\ul(\Rm)) \cap  
\CC^1(\Rm,L^2_\ul(\Rm))$ of \eqref{eq-u} such that, for all $L > 0$ 
and all $T > 0$, 
\[
  \sup_{t \in [-T,T]}\Bigl(
  \|u(x_n + \cdot,t_n + t) - \bar u(\cdot,t)\|_{H^1([-L,L])} +
  \|u_t(x_n + \cdot,t_n + t) - \bar u_t(\cdot,t)\|_{L^2([-L,L])}
  \Bigr) \,\xrightarrow[n \to \infty]{}\, 0~.
\]
\end{prop}

In other words, after extracting a subsequence, we can assume that the
sequence $\{u(x_n+x,t_n+t)\}$ converges in $\CC^0([-T,T],H^1_\loc(\Rm))
\cap \CC^1([-T,T],L^2_\loc(\Rm))$ towards a solution $\bar u(x,t)$ of
\eqref{eq-u}.

\begin{demo}
As in the proof of Proposition~\ref{propexist}, we set $w = (u,u_t)$ 
and we consider Eq.~\eqref{eq-w} instead of Eq.~\eqref{eq-u}. If $w_0 
\in X$, the solution of \eqref{eq-w} with initial data $w_0$ has the 
following representation:
\[
  w(t) \,=\, e^{At}w_0 + \int_0^t e^{A(t-s)}F(w(s))\dd s \,\equiv\,
  w_1(t) + w_2(t)~.
\]
As is easily verified, there exists $C_6 > 0$ and $\mu > 0$ such that
$\|e^{At}\|_{\LL(X)} \le C_6 \,e^{-\mu t}$ for all $t \ge 0$ (this
estimate will be established in a more general setting in
Section~\ref{sec9}, Lemma~\ref{Sinfty}). Thus $w_1(t) = e^{At}w_0$
converges exponentially to zero as $t \to +\infty$, and can therefore
be neglected.  On the other hand, by Proposition~\ref{propexist},
there exists $C_7 > 0$ such that $\|w(t)\|_X \le C_7$ for all $t \ge
0$. As $F$ maps $X$ into $Y= D(A)$ and is Lipschitz on bounded sets,
there exists $C_8 > 0$ such that $\|AF(w)\|_X \le C_8$ whenever
$\|w\|_X \le C_7$. Since $A w_2(t) = \int_0^t e^{A(t-s)} AF(w(s)) \dd
s$, we deduce that
\[
  \|Aw_2(t)\|_X \,\le\, C_6 \int_0^t e^{-\mu(t-s)} \|AF(w(s))\|_X
  \dd s \,\le\, \frac{C_6 C_8}{\mu}~, \quad t \ge 0~,
\]
hence there exists $C_9 > 0$ such that $\|w_2(t)\|_Y \le C_9$ for 
all $t \ge 0$. In particular, given any $T > 0$, the sequence
$\{w_2(x_n+\cdot,t_n-T)\}$ is bounded in $H^2([-L,L]) \times 
H^1([-L,L])$ for any $L > 0$. Extracting a subsequence and using 
a diagonal argument, we can assume that there exists $\bar w_0 \in 
H^2_\loc(\Rm) \times H^1_\loc(\Rm)$ such that, for any $L > 0$, 
\begin{equation}\label{w2loc}
  w_2(x_n+\cdot,t_n-T) \,\xrightarrow[n \to \infty]{}\, \bar w_0
  \quad \hbox{in} \quad H^1([-L,L]) \times L^2([-L,L])~.
\end{equation}
By construction $\|\bar w_0\|_Y \le C_9$, hence in particular $\bar
w_0 \in X$. Note that \eqref{w2loc} still holds if we replace $w_2$ by
the full solution $w$, because $\|w_1(\cdot,t)\|_X$ converges to zero.
Finally, let $\bar w(t) \in \CC^0([-T,+\infty),X)$ be the solution 
of \eqref{eq-w} with initial data $w(\cdot,-T) = \bar w_0$. Since 
the evolution defined by \eqref{eq-w} has a finite speed of 
propagation, it is clear that the solution $\bar w(t)$ depends 
continuously on the initial data $\bar w_0$ in the topology of
$H^1_\loc(\Rm) \times L^2_\loc(\Rm)$, uniformly in time on compact
intervals. Thus it follows from \eqref{w2loc} that, for all $L > 0$,
\[
  \sup_{t \in [-T,T]} \|w(x_n+\cdot,t_n+t) - \bar w(\cdot,t)
  \|_{H^1([-L,L]) \times L^2([-L,L])} \,\xrightarrow[n \to 
  \infty]{}\, 0~.
\]
Repeating the argument for larger $T$ and using a diagonal extraction, 
we conclude the proof of Proposition~\ref{propcompact}.
\end{demo}


\section{Pinching at infinity and splitting of the solution}
\label{sec3}

In this section we prove that, if the initial data satisfy the
boundary conditions \eqref{cond-},\eqref{cond+}, the solution $u(x,t)$
of Eq.~\eqref{eq-u} has the same properties for all positive times. As
a consequence, we show that any such solution can be decomposed into a
principal part $v(x,t)$ which is compactly supported to the right, and
a small remainder $r(x,t)$ which decays exponentially to zero as $t
\to +\infty$.

We first verify that, due to assumptions \eqref{u=0}, \eqref{u=1}, 
the homogeneous equilibria $u = 0$ and $u = 1$ are stable steady
states of Eq.~\eqref{eq-u}. Let $(u_0,u_1) \in X = H^1_\ul(\Rm) 
\times L^2_\ul(\Rm)$, and let $(u,u_t)$ be the solution of 
Eq.~\eqref{eq-u} with initial data $(u_0,u_1)$ given by 
Proposition~\ref{propexist}. 

\begin{lemma}\label{homostab}
There exist positive constants $K_i, \delta_i, \mu_i$ for $i = 0,1$
such that\\
{\bf a)} If $\|(u_0,u_1)\|_X^2 \le \delta_0$, then
$\|(u(\cdot,t),u_t(\cdot,t))\|_X^2 \le K_0 \|(u_0,u_1)\|_X^2 
\,e^{-\mu_0 t}$ for all $t \ge 0$.\\
{\bf b)} If $\|(u_0{-}1,u_1)\|_X^2 \le \delta_1$, then 
$\|(u(\cdot,t){-}1,u_t(\cdot,t))\|_X^2 \le K_1 \|(u_0{-}1,u_1)\|_X^2 
\,e^{-\mu_1 t}$ for all $t \ge 0$.
\end{lemma}

\begin{demo}
It is sufficient to prove {\bf a)}, the other case being similar. 
Let $\beta_0 = V''(0) > 0$, and choose $\epsilon_0 > 0$ small 
enough so that
\begin{equation}\label{beta0def}
  \frac{\beta_0}{2} \,\le\, V''(u) \,\le\, 2\beta_0~, \quad
  \hbox{for all } u \in [-\epsilon_0,\epsilon_0]~.
\end{equation}
In particular, we have
\begin{equation}\label{beta0deff}
  \frac{\beta_0 u^2}{2} \,\le\, u V'(u) \,\le\, 2\beta_0 u^2~, 
  \quad \hbox{and} \quad  
  \frac{\beta_0 u^2}{4} \,\le\, V(u) \,\le\, \beta_0 u^2~,
\end{equation}
whenever $|u| \le \epsilon_0$. In analogy with \eqref{EEdef}, we 
introduce the functional
\[
  \EE_0(\xi,t) \,=\, \int_\Rm (T_\xi \rho)(x)\Bigl(\alpha^2 u_t^2 
    + \alpha u_x^2 + 2\alpha V(u) + \frac12 u^2 
    + \alpha u u_t\Bigr)(x,t)\dd x~, 
\]
where $\rho(x) = \exp(-\kappa|x|)$ and $\kappa > 0$ is small enough 
so that $2\sqrt{\alpha}\kappa \le 1$ and $2\kappa^2 \le \beta_0$.
If $\|u(\cdot,t)\|_{L^\infty}^2 \le 2\|u(\cdot,t)\|_{H^1_\ul}^2 \le 
\epsilon_0^2$, it follows from \eqref{beta0deff} and from the 
definitions \eqref{L2ulnorm}--\eqref{Hkuldef} that
\[
   C_0^{-1} \|(u(\cdot,t),u_t(\cdot,t))\|_X^2 \,\le\, 
   \sup_{\xi\in\Rm} \EE_0(\xi,t) \,\le\, C_0 \|(u(\cdot,t),
   u_t(\cdot,t))\|_X^2~,
\]
for some $C_0 > 1$. Under the same assumption, we find as in 
the proof of Proposition~\ref{propexist}:
\begin{eqnarray*}
  \partial_t \EE_0(\xi,t) &=& -\int_\Rm 
    (T_\xi \rho)(\alpha u_t^2 + u_x^2 + u V'(u))\dd x 
    -\int_\Rm (T_\xi \rho)'(uu_x + 2\alpha u_x u_t)\dd x \\
  &\le& -\int_\Rm (T_\xi \rho)(\frac\alpha2 u_t^2 + 
    \frac12 u_x^2 + \frac{\beta_0}4 u^2)\dd x \,\le\, 
  -\mu_0 \EE_0(\xi,t)~,
\end{eqnarray*}
for some $\mu_0 > 0$. Now, let $K_0 = C_0^2$ and choose $\delta_0 > 0$
small enough so that $2 K_0 \delta_0 < \epsilon_0^2$. If
$\|(u_0,u_1)\|_X^2 \le \delta_0$, the inequalities above imply that
the solution $(u,u_t)$ satisfies $\|(u(\cdot,t), u_t(\cdot,t))\|_X^2
\le K_0 \|(u_0,u_1) \|_X^2 \,e^{-\mu_0 t}$ for all $t \ge 0$.  In
particular, $\|u(\cdot,t)\|_{L^\infty}^2 \le
2\|u(\cdot,t)\|_{H^1_\ul}^2 \le \epsilon_0^2 \,e^{-\mu_0 t}$ for all
$t \ge 0$.
\end{demo}

From now on, we assume that the initial data $(u_0,u_1) \in X$
satisfy the assumptions \eqref{cond-}, \eqref{cond+} for some 
$\delta \le \min(\delta_0,\delta_1)/2$, and we let $u \in 
\CC^0([0,+\infty),H^1_\ul(\Rm)) \cap \CC^1([0,+\infty),L^2_\ul(\Rm))$
be the solution of \eqref{eq-u} given by Theorem~\ref{mainthm}. 
Using Lemma~\ref{homostab} and the finite speed of propagation 
we show that, for all $t \ge 0$, the solution $u(x,t)$ stays close 
for large $|x|$ to the homogenous equilibria $u = 0$ and $u = 1$. 

\begin{prop}\label{pinching}
If $\delta \le \min(\delta_0,\delta_1)/2$, the solution of 
\eqref{eq-u} given by Theorem~\ref{mainthm} satisfies, for 
all $t \ge 0$, 
\begin{eqnarray}\label{pinch+}
  \limsup_{\xi \to +\infty} \int_\xi^{\xi+1} \Bigl(u(x,t)^2 + {u_x(x,t)}^2
  + u_t(x,t)^2\Bigr) \dd x &\le& K_0 \delta_0\,e^{-\mu_0 t}~,
  \\ \label{pinch-}
  \limsup_{\xi \to -\infty} \int_\xi^{\xi+1} \Bigl((u(x,t)-1)^2 + 
  {u_x(x,t)}^2 + u_t(x,t)^2\Bigr) \dd x &\le& K_1 \delta_1\,e^{-\mu_1 t}.
\end{eqnarray}
\end{prop}

\begin{demo}
We only prove the first inequality, the second one being similar. 
Take $\xi_0 \in \Rm$ such that
\[
  \int_\xi^{\xi+1} \Bigl(u_0(x)^2 + {u_0'(x)}^2
  + u_1(x)^2\Bigr) \dd x \,\le\,  \frac{3\delta_0}{4}~,
  \quad \hbox{for all } \xi \,\ge\, \xi_0-4~.
\]
We consider the modified initial data $(r_0,r_1) \in X$ defined by 
\begin{equation}\label{r0r1def}
  r_0(x) \,=\, \theta(x-\xi_0)u_0(x)~, \quad
  r_1(x) \,=\, \theta(x-\xi_0)u_1(x)~, \quad
  x \in \Rm~,
\end{equation}
where $\theta(x) = \min(1,(1+x/4)_+)$ satisfies $\theta(x) = 1$ 
for $x \ge 0$, $\theta(x) = 0$ for $x \le -4$, and $|\theta'(x)| 
\le 1/4$ for all $x$. By construction $(r_0(x),r_1(x)) = 
(u_0(x),u_1(x))$ for all $x \ge \xi_0$, and
\[
  \|(r_0,r_1)\|_X^2 \,\le\, \frac43 \sup_{\xi \ge \xi_0-4}
  \int_\xi^{\xi+1} \Bigl(u_0(x)^2 + {u_0'(x)}^2
  + u_1(x)^2\Bigr) \dd x \,\le\, \delta_0~.
\]
If $(r,r_t) \in \CC^0([0,+\infty),X)$ is the solution of \eqref{eq-u}
with initial data $(r_0,r_1)$, we know from Lemma~\ref{homostab}
that 
\begin{equation}\label{rbound1}
  \|(r(\cdot,t),r_t(\cdot,t))\|_X^2 \,\le\, K_0 \delta_0 \,e^{-\mu_0 t}~,
  \quad \hbox{for all } t \ge 0~.
\end{equation}
On the other hand, the finite speed of propagation implies that
$u(x,t) = r(x,t)$ for all $t \ge 0$ and all $x \ge \xi_0 +
t/\sqrt{\alpha}$. Both observations together imply \eqref{pinch+}.
\end{demo}

\medskip\noindent {\bf Decomposition of the solution:} 
The proof of Proposition~\ref{pinching} provides us with a 
useful decomposition of the solution of \eqref{eq-u}. Let
\begin{equation}\label{udecomp}
  u(x,t) \,=\, v(x,t) + r(x,t)~, \quad x \in \Rm~, \quad t \ge 0~,
\end{equation}
where $r(x,t)$ is the solution of \eqref{eq-u} associated to 
the initial data $(r_0,r_1)$ defined in \eqref{r0r1def}. By 
construction, the principal part $v(x,t)$ vanishes identically 
for $x \ge \xi_0 + t/\sqrt{\alpha}$, and satisfies the modified 
equation
\begin{equation}\label{eq-v}
   \alpha v_{tt} + v_t \,=\, v_{xx} -V'(v+r) + V'(r)~,
\end{equation}
supplemented with the initial data $(v_0,v_1) = (u_0-r_0,u_1-r_1)$. 
If we define
\begin{equation}\label{fdef}
   f(v,r) \,=\, V'(v) + V'(r) - V'(v+r) \,=\, 
   -vr \int_0^1\!\int_0^1 V'''(tv+sr)\dd t\dd s~,
\end{equation}
we can rewrite \eqref{eq-v} in the form
\begin{equation}\label{eq-v2}
   \alpha v_{tt} + v_t \,=\, v_{xx} - V'(v) + f(v,r)~.
\end{equation}
The main advantage of working with \eqref{eq-v2} instead of
\eqref{eq-u} is that the energy functional \eqref{Ecdef} (with $u$
replaced by $v$) is well-defined for all $c > 0$ since $v(x,t)$ is
compactly supported to the right. The price to pay is the additional
term $f(v,r)$ in \eqref{eq-v2}, which we shall treat as a perturbation.
Remark that, since $v(x,t)$ and $r(x,t)$ stay uniformly bounded for
all $t \ge 0$, the formula \eqref{fdef} shows that there exists $K_2 >
0$ such that
\begin{equation}\label{fbound}
  |f(v(x,t),r(x,t))| \,\le\,  K_2 \,|v(x,t)| \,|r(x,t)|~, 
  \quad x \in \Rm~, \quad t \ge 0~.
\end{equation}
Moreover, we know that $\|r(\cdot,t)\|_{L^\infty}^2 \le 2\|r(\cdot,t)
\|_{H^1_\ul}^2 \le \epsilon_0^2 \,e^{-\mu_0 t}$ for all $t \ge 0$,
hence \eqref{eq-v2} is really a small perturbation of \eqref{eq-u} 
for large times. In particular, the asymptotic compactness 
property stated in Proposition~\ref{propcompact} holds for the
solution $v(x,t)$ of \eqref{eq-v2}, and by Proposition~\ref{propexist}
there exists $M_0 > 0$ such that 
\begin{equation}\label{vbound}
  \|v(\cdot,t)\|_{H^1_\ul}^2 + \|v_t(\cdot,t)\|_{L^2_\ul}^2 
  \,\le\, M_0^2~, \quad \hbox{for all }t \ge 0~. 
\end{equation}

\medskip\noindent {\bf The invasion point:} 
As is explained in \cite{GR5,Ri}, to control the behavior of the
solution $u(x,t)$ of \eqref{eq-u} using the energy functionals
\eqref{Ecdef} it is necessary to track for all times the approximate
position of the front interface. Since $r(x,t)$ converges uniformly to
zero as $t \to +\infty$, this can be done for the solution $v(x,t)$ 
of \eqref{eq-v2} instead of $u(x,t)$. We thus introduce the 
invasion point $\overline x(t) \in \Rm$ defined for any $t \ge 0$ by
\begin{equation}\label{barxdef}
  \overline x(t) \,=\, \sup\{x \in \Rm \,|\, |v(x,t)| \ge \epsilon_0\}~,
\end{equation}
where $\epsilon_0$ is as in \eqref{beta0def}. It is clear that
$\overline x(t) < \xi_0 + t/\sqrt{\alpha}$ since $v(x,t)$ vanishes
identically for larger values of $x$. In the same way, using
\eqref{pinch-}, one can prove that there exists $\xi_1 \in \Rm$ such
that $\overline x(t) > \xi_1 - t/\sqrt{\alpha}$ for all $t \ge 0$.
Note that $\overline x(t)$ is not necessarily a continuous function of
$t$, although it follows from the definition \eqref{barxdef} that
$\overline x(t)$ is upper semi-continuous.


\section{Energy estimates in a Galilean frame}
\label{sec4}

As is explained in the introduction, the proof of 
Theorem~\ref{mainthm} is based on the existence of Lyapunov
functions for Eq.~\eqref{eq-u} in uniformly translating frames. 
The aim of this section is to define these functions rigorously 
and to study their basic properties. 

Let $u(x,t)$ be a solution of \eqref{eq-u} whose initial data satisfy
the assumptions of Theorem~\ref{mainthm}. Given any $c > 0$, we go
to a uniformly translating frame by setting, as in \eqref{def-uc},
$u(x,t) = u_c(\sqrt{1+\alpha c^2}\,x - ct\,,\,t)$. To avoid
confusions, we always denote by $y = \sqrt{1+\alpha c^2}\,x - ct$ 
the space variable in the moving frame. Note that the 
{\em physical speed} $s \in (0,1/\sqrt{\alpha})$ of the 
frame is related to the {\em parabolic speed} $c > 0$ by the
formulas
\begin{equation}\label{eq-s-c}
  s \,=\, \frac c{\sqrt{1+\alpha c^2}}~, \qquad
  c \,=\, \frac s{\sqrt{1-\alpha s^2}}~.
\end{equation}
If $u(x,t)$ is decomposed according to \eqref{udecomp}, then 
$u_c(y,t) = v_c(y,t) + r_c(y,t)$ where
\begin{equation}\label{def-vc-rc}
  v_c(y,t) \,=\, v\Bigl(\frac{y+ct}{\sqrt{1+\alpha c^2}}\,,\,t
  \Bigr)~, \quad \hbox{and} \quad 
  r_c(y,t) \,=\, r\Bigl(\frac{y+ct}{\sqrt{1+\alpha c^2}}\,,\,t \Bigr)~.
\end{equation}
By construction, both $v_c$ and $r_c$ belong to $\CC^0([0,+\infty),
H^1_\ul(\Rm)) \cap \CC^1([0,+\infty),L^2_\ul(\Rm))$. Moreover,
from \eqref{vbound} and Lemma~\ref{homostab}, we know that
\begin{equation}\label{vcrcbound}
  \|v_c(t)\|_{L^\infty}^2 \,\le\, 2\|v(t)\|_{H^1_\ul}^2 
  \,\le\, M_0^2~, \quad 
  \|r_c(t)\|_{L^\infty}^2 \,\le\, 2\|r(t)\|_{H^1_\ul}^2 
  \,\le\, \epsilon_0^2 \,e^{-\mu_0 t}~,
\end{equation}
for all $t \ge 0$. In view of \eqref{eq-uc}, \eqref{eq-v2}, the 
evolution equations satisfied by $v_c$, $r_c$ read
\begin{equation}\label{eq-vc-rc}
  \left\{\begin{array}{l} 
  \alpha \ddot r_c + \dot r_c - 2\alpha c \dot r'_c \,=\, r''_c + 
  c r'_c-V'(r_c)~,\\
  \alpha \ddot v_c + \dot v_c - 2\alpha c \dot v'_c \,=\, v''_c
  +c v'_c - V'(v_c) + f(v_c,r_c)~,
  \end{array}\right.
\end{equation}
where $f(v_c,r_c) = -V'(v_c + r_c)+V'(v_c)+V'(r_c)$. Here and in the
rest of the text, to simplify the notation and to avoid double 
subscripts, we denote $\dot v_c(y,t) = \partial_t v_c(y,t)$, 
$v_c'(y,t) = \partial_y v_c(y,t)$, and similarly for $r_c$. 
In analogy with \eqref{barxdef}, we also define the {\em invasion
point} in the moving frame by 
\begin{equation}\label{barydef}
  \overline y_c(t) \,=\, \sqrt{1+\alpha c^2}\,\overline x(t)-ct 
  \,=\, \sup\{y \in \Rm \,|\, |v_c(y,t)| \ge \epsilon_0\}~.
\end{equation}

\subsection{The energy functional}

In a moving frame with parabolic speed $c > 0$, the energy functional
involves an exponentially growing weight $e^{cy}$, see \eqref{Ecdef}. 
It is thus natural to introduce the following weighted spaces:
\begin{align}\label{def-espaces}
  L^2_c(\Rm) &\,=\, \{u \in L^2_\loc(\Rm) \,|\, e^{cy/2} u \in 
    L^2(\Rm)\}~, \\  \nonumber
  H^1_c(\Rm) &\,=\, \{u \in H^1_\loc(\Rm)\,|\, e^{cy/2} u \in 
    L^2(\Rm) \hbox{ and } e^{cy/2}u' \in L^2(\Rm)\}~.
\end{align} 
Since $v_c(\cdot,t) \in H^1_\ul(\Rm)$ and $v_c(y,t)$ vanishes for all
sufficiently large $y > 0$, it is clear that $v_c(\cdot,t)$ belongs to
$H^1_c(\Rm)$ for any $c > 0$. Similarly, $\dot v_c(\cdot,t)$ belongs
to $L^2_c(\Rm)$. The following quantity is thus well-defined for 
any $y_0 \in \Rm$ and all $t \ge 0$:
\begin{equation}\label{def-energy}
  E_c(y_0,t) \,=\, \int_\Rm e^{cy}\Bigl(\frac\alpha2 |\dot v_c|^2 +
  \frac12 |v'_c|^2 + V(v_c) \Bigr)(y_0+y,t)\dd y~.
\end{equation}
We shall refer to $E_c(y_0,t)$ as the {\em energy} of the solution
$v_c(y,t)$ in the moving frame. The translation parameter $y_0$ is
introduced here for later convenience. Changing $y_0$ results in a
simple rescaling, as is clear from the identity
\begin{equation}\label{eq-Ec-y0}
  E_c(y_0,t) \,=\, e^{c(y_1-y_0)}E_c(y_1,t)~.
\end{equation}
Due to the term $f(v_c,r_c)$ in \eqref{eq-vc-rc}, the energy 
$E_c(y_0,t)$ is not necessarily a decreasing function of time. 
Indeed, a formal calculation gives
\begin{equation}\label{energy-decay}
  \partial_t E_c(y_0,t) \,=\, -(1+\alpha c^2)\int_\Rm e^{cy} |\dot
  v_c(y_0+y,t)|^2 \dd y + R_c(y_0,t)~,
\end{equation}
where 
\begin{equation}\label{Rcdef}
  R_c(y_0,t) \,=\, \int_\Rm e^{cy} (f(v_c,r_c) \dot v_c)(y_0+y,t)
  \dd y~.
\end{equation}
Using \eqref{fbound} and \eqref{vcrcbound}, it is easy to verify 
that $R_c(y_0,t)$ is well-defined for all $t \ge 0$ and 
depends continuously on time. A classical argument then shows
that $E_c(y_0,t)$ is indeed differentiable with respect to 
$t$ and that \eqref{energy-decay} holds for all $t \ge 0$. 
The purpose of this section is to show that, in appropriate
situations, the remainder term $R_c$ in \eqref{energy-decay} is a
negligible quantity which does not really affect the decay of the
energy. Our first result in this direction is:

\begin{lemma}\label{lemme-estim-reste}
There exists a positive constant $K_3$, independent of
$c$, such that
\begin{equation}\label{Rcest1}
  |R_c(y_0,t)| \,\le\, K_3 \,e^{-\mu t} \Bigl(E_c(y_0,t) + 
  \frac{K_3}c \,e^{c(\overline y_c(t)-y_0)}\Bigr)~,
\end{equation}
for all $y_0 \in \Rm$ and all $t \ge 0$, where $\mu = \mu_0/2$. 
\end{lemma}

\begin{demo}
Using \eqref{fbound}, \eqref{vcrcbound}, and \eqref{Rcdef}, we obtain
\begin{align*}
  |R_c(y_0,t)| &\,\le\, K_2 \epsilon_0 \,e^{-\mu t} \!\int_\Rm 
    e^{c(y-y_0)}|\dot v_c \,v_c|(y,t)\dd y \\ 
  &\,\le\, \frac{K_2\epsilon_0}2 \,e^{-\mu t} \!\int_\Rm e^{c(y-y_0)}
  (|\dot v_c|^2+|v_c|^2)(y,t)\dd y~.
\end{align*}
If $y \ge \overline y_c(t)$, then $|v_c(y,t)| \le \epsilon_0$ by
\eqref{barydef}, hence $|v_c(y,t)|^2 \le (4/\beta_0) V(v_c(y,t))$ by
\eqref{beta0deff}. Thus
\[
  \frac12 (|\dot v_c(y,t)|^2+|v_c(y,t)|^2) \,\le\, 
  C \Bigl(\frac{\alpha}2|\dot v_c(y,t)|^2 + \frac12 |v'_c(y,t)|^2
  +V(v_c(y,t))\Bigr)~,
\]
where $C = \max(\alpha^{-1},2\beta_0^{-1})$. If $y \le \overline 
y_c(t)$, we can bound
\begin{align*}
  \frac12 (|\dot v_c(y,t)|^2+|v_c(y,t)|^2) \,\le\, 
  &C \Bigl(\frac{\alpha}2|\dot v_c(y,t)|^2 + \frac12 |v'_c(y,t)|^2
  +V(v_c(y,t))\Bigr)\\ 
  &+ \frac12 \|v_c(t)\|_{L^\infty}^2 + C |\min V|~.
\end{align*}
Combining these estimates and using \eqref{vcrcbound}, we thus 
obtain
\begin{align*}
  |R_c(y_0,t)| &\,\le\, K_2 \epsilon_0 \,e^{-\mu t} \left(
  C E_c(y_0,t) + \Bigl(\frac{M_0^2}2 + C |\min V|\Bigr)
  \int_{-\infty}^{\overline y_c(t)} e^{c(y-y_0)}\dd y\right) \\
  &\,\le\, K_3 \,e^{-\mu t}\left(E_c(y_0,t)+\frac{K_3}{c} 
  \,e^{c(\overline y_c(t)-y_0)}\right)~,
\end{align*}
which is the desired result.
\end{demo}

The following corollary of Lemma~\ref{lemme-estim-reste} will turn 
out to be useful: 

\begin{prop}\label{prop-estim-reste}
Assume that the invasion point satisfies, for some $c_+ > 0$, 
\[
  \limsup_{t\to +\infty} \frac {\overline y_{c_+}(t)}t \,\le\, 0~.
\]
Then, there exist $\eta > 0$ and $K_4 > 0$ such that, for all
$c \in [c_+-\eta,c_++\eta]$, all $y_0 \in \Rm$, and all $t_1 \ge 
t_0 \ge 0$,
\begin{equation}\label{Ecbd1}
  E_c(y_0,t_1) \,\le\, K_4 \max(E_c(y_0,t_0)\,,\,e^{-cy_0})~.
\end{equation}
Moreover
\begin{equation}\label{Rcest2}
  |R_c(y_0,t)| \,\le\, K_4 \,e^{-\mu t/2} \max(E_c(y_0,t_0)\,,\,
  e^{-cy_0})~, \quad \hbox{for all } t \ge t_0~.
\end{equation}
\end{prop}

\begin{demo}
By assumption, for all $\epsilon>0$, there exists $C_\epsilon > 0$
such that $\overline y_{c_+}(t) \le \epsilon t + C_\epsilon$
for all $t \ge 0$. In a frame moving at parabolic speed $c$, 
this bound becomes
\[
  \overline y_c(t) \,\le\, \left(\sqrt{\frac{1+\alpha c^2}{1+\alpha
      c^2_+}} (c_++\epsilon)-c\right)t ~+~ \sqrt{\frac{1+\alpha
    c^2}{1+\alpha c^2_+}} \,C_\epsilon~.
\]
Thus, if we choose $\epsilon > 0$ small enough, there exist $\eta 
\in (0,c_+)$ and $C_1 > 0$ such that, for all $c \in [c_+-\eta,
c_++\eta]$ and all $t \ge 0$, we have $c\overline y_c(t) \le 
(\mu t)/2 + C_1$. Using \eqref{energy-decay} and 
Lemma~\ref{lemme-estim-reste}, we deduce that
\begin{equation}\label{dem-prop-estim-reste-2}
  \partial_t E_c(y_0,t) \,\le\, |R_c(y_0,t)| \,\le\, 
  K_3 \,e^{-\mu t}(E_c(y_0,t)+ C_2 \,e^{\mu t/2 - cy_0})~,
\end{equation}
for some $C_2 > 0$. Integrating this differential inequality
between $t_0$ and $t_1$, we obtain
\begin{align*}
  E_c(y_0,t_1) &\,\le\, e^{\frac{K_3}{\mu}(e^{-\mu t_0}- e^{-\mu t_1})} 
  E_c(y_0,t_0) + K_3 C_2 \int_{t_0}^{t_1} e^{\frac{K_3}{\mu}(e^{-\mu t}
   - e^{-\mu t_1})} \,e^{-\mu t/2 -cy_0}\dd t \\
  &\,\le\, K_4 \max(E_c(y_0,t_0)\,,\,e^{-cy_0})~,
\end{align*}
which proves \eqref{Ecbd1}. Estimate \eqref{Rcest2} is a direct 
consequence of \eqref{Ecbd1} and \eqref{dem-prop-estim-reste-2}. 
\end{demo}

\begin{rem}
Of course, if the initial data $u_0,u_1$ decay rapidly enough as $x
\to +\infty$, the decomposition \eqref{udecomp} is not needed and we can
use the energy \eqref{Ecdef} instead of \eqref{def-energy}. This is the
point of view adopted in \cite{GR5}. In a first reading of the paper,
it is therefore possible to set $r = 0$ everywhere, in which case the
remainder term $R_c(y_0,t)$ disappears from \eqref{energy-decay} and
the energy $E_c(y_0,t)$ is a true Lyapunov function. However, once 
the the invasion point is under control, the results of this section 
show that $R_c(y_0,t)$ becomes really negligible for large times. 
Thus the general case can be seen as a perturbation of the particular
situation where $r = 0$, and the outline of the arguments is the
same in both cases.
\end{rem}

\subsection{A Poincar\'e inequality}

As was already noted in \cite{Muratov} and \cite{GR5}, Poincar\'e
inequalities hold in the weighted space $H^1_c(\Rm)$ if $c > 0$. 

\begin{prop}\label{prop-Poincare}
Let $c > 0$ and let $v_c \in H^1_c(\Rm)$. Then $e^{cy}|v_c(y)|^2 \to 0$ 
as $y \to +\infty$. Moreover, for any $y_1 \in \Rm$,
\begin{equation}\label{eq-Poincare}
  \frac{c^2}4 \int_{y_1}^\infty e^{cy}|v_c(y)|^2 \dd y \,\le\, 
  \int_{y_1}^\infty e^{cy} |v'_c(y)|^2\dd y~,
\end{equation}
and 
\begin{equation}\label{eq-Poincare-2}
  ce^{cy_1}|v_c(y_1)|^2 \,\le\,\int_{y_1}^\infty e^{cy}|v'_c(y)|^2 
  \dd y~.
\end{equation}
\end{prop}

\begin{demo}
A simple integration shows that, for all $y_1 \le y_2$, 
\[
  e^{cy_2}|v_c(y_2)|^2-e^{cy_1}|v_c(y_1)|^2 \,=\, 2 \int_{y_1}^{y_2}
  e^{cy}v'_c(y) v_c(y)\dd y + c \int_{y_1}^{y_2} e^{cy}|v_c(y)|^2 
  \dd y~.
\]
When $y_2$ goes to $+\infty$, both integrals in the right-hand side
have a finite limit since $v_c \in H^1_c(\Rm)$. Thus the first 
term in the left-hand side also has a limit, which is necessarily 
zero since $y \mapsto e^{cy}|v_c(y)|^2 \in L^1(\Rm)$. It follows that
\begin{equation}\label{dem-eq-poincare}
  e^{cy_1}|v_c(y_1)|^2 \,\le\,  2 \int_{y_1}^\infty e^{cy}|v'_c(y) 
  v_c(y)|\dd y - c \int_{y_1}^\infty e^{cy}|v_c(y)|^2 \dd y~.
\end{equation}
Now, for any $d > -c$, we have $|2v_c v'_c| \le (c+d)|v_c|^2 + 
\frac{1}{c+d}|v'_c|^2$. Inserting this bound in
\eqref{dem-eq-poincare} we find
\[
  e^{cy_1}|v_c(y_1)|^2\leq \frac 1{c+d}\int_{y_1}^\infty
  e^{cy}|v'_c(y)|^2\dd y+ d 
  \int_{y_1}^\infty e^{cy}|v_c(y)|^2 \dd y~,
\]
from which \eqref{eq-Poincare} follows by taking $d = -c/2$ 
and \eqref{eq-Poincare-2} by choosing $d = 0$.
\end{demo}

The Poincar\'e inequality implies the following important lower 
bound on the energy. We recall that, for all $y \ge \overline 
y_c(t)$, one has $|v_c(y,t)| \le \epsilon_0$ by \eqref{barydef}, 
so that $V(v_c(y,t)) \ge 0$ by \eqref{beta0deff}. Thus
\begin{align*}
  E_c(y_0,t) &\,=\, e^{-cy_0}\int_{-\infty}^{\overline y_c(t)}
  e^{cy}\Bigl(\frac\alpha2 |\dot v_c|^2 + \frac12 |v'_c|^2 + 
  V(v_c)\Bigr)(y,t)\dd y\\
  &\quad + e^{-cy_0}\int^{\infty}_{\overline y_c(t)} e^{cy}\Bigl
  (\frac\alpha2 |\dot v_c|^2 + \frac12 |v'_c|^2 + V(v_c) \Bigr)
  (y,t)\dd y\\
  &\,\ge\, e^{-cy_0} \int_{-\infty}^{\overline y_c(t)} e^{cy} 
  (\min V)\dd y + e^{-cy_0}\frac12 \int^{\infty}_{\overline y_c(t)}
  e^{cy} |v'_c(y,t)|^2 \dd y~.
\end{align*}
Using now \eqref{eq-Poincare-2} and the fact that $|v_c(\overline 
y_c(t),t)| = \epsilon_0$, we obtain
\begin{equation}\label{eq-lowerbound}
  E_c(y_0,t) \,\ge\, e^{c(\overline y_c(t)-y_0)}\Bigl(\frac 
  {\min V}c + \frac{c\epsilon_0^2}2 \Bigr)~. 
\end{equation}


\section{Existence of the invasion speed}
\label{sec5}

The purpose of this section is to show that the invasion point
$\overline x(t)$ defined in \reff{barxdef} has a positive 
average speed as $t \to +\infty$:

\begin{prop}\label{prop-asympt-speed}
The limit $\displaystyle \lim_{t \to +\infty} \frac{\overline
x(t)}{t} = s_\infty$ exists and lies in the interval
$(0,\frac1{\sqrt{\alpha}})$.
\end{prop}

\noindent
We call $s_\infty$ the {\em invasion speed} because this is 
the speed at which the front interface described by the 
solution $u(x,t)$ ``invades'' the steady state $u = 0$. 
We prove that $s_\infty < \frac 1{\sqrt{\alpha}}$, which means 
that the invasion process is always ``subsonic''. This 
might not be the case if we drop the assumption that the 
equilibrium $u = 0$ is stable. For instance, if $h(x) = 
(1+e^x)^{-1}$, one can check that $u(x,t) = h(x-st)$ is a solution 
of \eqref{eq-u} provided that 
\[ 
  -V'(u) \,=\, u(1-u)(s+\gamma(1-2u))~, \quad \hbox{where }
  \gamma \,=\,\alpha s^2-1~.
\]
If we choose $s > 0$ large enough so that $\gamma > 0$, the 
front $h(x-st)$ is supersonic, but in that case we also have
$V''(0) < 0$, hence $u = 0$ is an unstable equilibrium of 
\eqref{eq-u}. 

Our proof of Proposition~\ref{prop-asympt-speed} follows closely 
the method introduced in \cite{Ri} and simplified in \cite{GR5}. 
It is divided into three lemmas. 

\begin{lemma}
One has $~\displaystyle \limsup_{t \to +\infty} \frac{\overline x(t)}t
< \frac 1{\sqrt{\alpha}}~$.
\end{lemma}

\begin{demo}
Choose $c > 0$ large enough so that $\frac{\min V}c + 
\frac{c\epsilon_0^2}2 > 0$. By \eqref{eq-lowerbound}, 
there exists $C_1 > 0$ such that $E_c(0,t) \ge C_1\,e^{c\overline 
y_c(t)}$ for all $t \ge 0$. Inserting this bound into 
\eqref{Rcest1}, we see that there exists $C_2 > 0$ such that
\[
  \partial_t E_c(0,t) \,\le\, |R_c(0,t)| \,\le\, C_2 
  \,e^{-\mu t}E_c(0,t)~, \quad \hbox{for all } t \ge 0~.
\]
If we integrate this differential inequality as in the proof of 
Proposition \ref{prop-estim-reste}, we find that $E_c(0,t) 
\le C_3 E_c(0,0)$ for all $t \ge 0$. Going back to the lower
bound $E_c(0,t) \ge C_1\,e^{c\overline y_c(t)}$, we conclude
that $\overline y_c(t)$ is bounded from above, hence
\[
  \limsup_{t\to +\infty}\frac{\overline x(t)}t \,=\, 
  \frac 1{\sqrt{1+\alpha c^2}}\Bigl(c + \limsup_{t\to +\infty}
  \frac{\overline y_c(t)}t\Bigr) \,\le\, \frac c{\sqrt{1+\alpha c^2}}
  \,<\, \frac1{\sqrt{\alpha}}~,
\]
which is the desired result.
\end{demo}

\begin{lemma}\label{lemme-limsup-positive}
One has $~\displaystyle \limsup_{t\rightarrow+\infty} \frac{\overline 
x(t)}t > 0~$.
\end{lemma}

\begin{demo}
We argue by contradiction and assume that $\limsup(\overline y_c(t)/t)
< 0$ for all $c > 0$. Using \eqref{energy-decay}, \eqref{Rcdef}
together with the bound \eqref{fbound}, we find 
\begin{align*}
  \partial_t E_c(0,t) &\,=\, -(1+\alpha c^2)\int_{\Rm} e^{cy}|\dot 
  v_c|^2(y,t)\dd y + \int_\Rm e^{cy} (f(v_c,r_c)\dot v_c)(y,t)\dd y\\
  &\,\le\, \frac14 \int_{\Rm} e^{cy} |f(v_c,r_c)|^2(y,t)\dd y 
  \,\le\, \frac{K_2^2}4 \int_\Rm e^{cy} |v_c(y,t)|^2 |r_c(y,t)|^2\dd y~.
\end{align*}
Our goal is to bound the right-hand side by a quantity which is
integrable in time and independent of $c$ if $c$ is sufficiently 
small. To do that, we fix $c' = 2\sqrt{\mu}$, where $\mu > 0$ is
as in Lemma~\ref{lemme-estim-reste}, and we assume that 
$c \in (0,c']$. Denoting $\rho = \sqrt{1+\alpha c^2}/\sqrt{1+\alpha 
{c'}^2}$ and using the definitions \eqref{def-vc-rc}, we obtain
the identity
\begin{align*}
  \int_\Rm e^{cy} |v_c|^2 |r_c|^2(y,t) \dd y \,&=\, \int_\Rm
  e^{cy} |v_{c'}|^2 |r_{c'}|^2 (\rho^{-1}(y+ct)-c't\,,\,t)\dd y \\ 
  &\,=\, \rho \,e^{c(\rho c'-c)t} \int_\Rm e^{c\rho z}
  |v_{c'}|^2|r_{c'}|^2(z,t) \dd z~.
\end{align*}
Since $\rho \le 1$ and $c(\rho c'-c) \le \mu$ by construction, 
we have
\[
  \int_\Rm e^{cy}|v_c(y,t)|^2|r_c(y,t)|^2\dd y \,\le\, e^{\mu t}
  \left(\int_0^\infty e^{c'y}|v_{c'}|^2|r_{c'}|^2(y,t)\dd y + 
  \int_{-\infty}^0 |v_{c'}|^2|r_{c'}|^2(y,t)\dd y \right)~.
\]
Remark that the right-hand side is now independent of $c$. To 
bound the first integral, we proceed as in the proof of 
Lemma~\ref{lemme-estim-reste}. Since $\overline y_{c'}(t)$ 
is bounded from above by assumption, so is $E_{c'}(0,t)$ 
by Proposition~\ref{prop-estim-reste} and we obtain
\begin{align*}
  &\int_0^\infty e^{c'y}|v_{c'}|^2|r_{c'}|^2(y,t)\dd y \,\le\, 
  \epsilon_0^2\,e^{-2\mu t} \int_\Rm e^{c'y}|v_{c'}(y,t)|^2\dd y\\
  &\qquad \,\le\, C \epsilon_0^2\,e^{-2\mu t} \Bigl(E_{c'}(0,t) + 
  (M_0^2 + |\min V|) \frac1{c'}\, e^{c'\overline y_{c'}(t)}\Bigr)
  \,\le\, C_4 \,e^{-2\mu t}~,
\end{align*}
for some $C_4 > 0$. To estimate the second integral we observe that 
$r(x,t) = 0$ for $x \le \xi_0-4-t/\sqrt{\alpha}$, because the 
initial data $(r_0,r_1)$ satisfy \eqref{r0r1def}. Thus there 
exists $C_5 > 0$ such that $r_{c'}(y,t) = 0$ whenever $y 
\le -C_5(1+t)$, hence
\[
  \int_{-\infty}^0 |v_{c'}|^2|r_{c'}|^2 (y,t)\dd y \,\le\, 
  C_5(1+t) \|v_{c'}(t)\|^2_{L^\infty} \|r_{c'}(t)\|^2_{L^\infty}
  \,\le\, C_5(1+t) M_0^2 \epsilon_0^2 \,e^{-2\mu t}~.
\] 
Summarizing, we have shown the existence of a constant $C_6 > 0$ such 
that $\partial_t E_c(0,t) \le C_6(1+t)\,e^{-\mu t}$
for all $t \ge 0$ and all $c \in (0,c']$. In particular, 
$\int_0^\infty \partial_t E_c(0,t)\dd t \le C_7 = C_6(1+\mu)/
\mu^2$. 

Now, if the initial data $(u_0,u_1)$, or equivalently $(v_0,v_1)$, 
satisfy the boundary condition \eqref{cond-} for some sufficiently
small $\delta > 0$, it is straightforward to verify that 
\begin{equation}\label{invasion_cond}
  \frac{E_c(0,0)}{\sqrt{1+\alpha c^2}} \,=\, \int_\Rm 
  e^{c\sqrt{1+\alpha c^2}\,x} \Bigl(\frac\alpha2 |v_1{+}sv_0'|^2
  + \frac{1}{2(1{+}\alpha c^2)} |v_0'|^2 + V(v_0)\Bigr)(x)\dd x
  \,\xrightarrow[c\to 0]{}\, -\infty~, 
\end{equation}
because $V(1) < 0$. In particular we can take $c \in (0,c']$
small enough so that $E_c(0,0) \le -2C_7$. Then $E_c(0,t) \le
-C_7$ for all $t \ge 0$, and since $E_c(0,t) \ge -e^{c\overline
y_c(t)}|\min V|/c$ by \eqref{eq-lowerbound}, we conclude that 
$\overline y_c(t)$ is bounded from below. Thus $\limsup 
(\overline y_c(t)/t) \ge 0$, which is the desired contradiction. 
\end{demo}

\begin{lemma}\label{lemme-asympt-speed}
One has $~\displaystyle \liminf_{t\rightarrow+\infty} \frac{\overline 
x(t)}t = \limsup_{t\rightarrow+\infty} \frac{\overline x(t)}t~$.
\end{lemma}

\begin{demo}
Again, we argue by contradiction and assume that
\[
   s_- \,\equiv\, \liminf_{t \to +\infty}\frac{\overline x(t)}t 
   \,<\, \limsup_{t \to +\infty} \frac{\overline x(t)}t 
   \,\equiv\, s_+~.
\]
Then there exist two increasing sequences of times $t_n$ and $t'_n$, 
both converging to $+\infty$, such that
\[
  \frac{\overline x(t_n)}{t_n}\,\xrightarrow[n\to\infty]{}\,s_+~, 
  \quad \hbox{and} \quad \frac{\overline x(t'_n)}{t'_n} 
  \,\xrightarrow[n\to\infty]{}\,s_-~. 
\]
Given $T > 0$, we can assume in view of Proposition \ref{propcompact}
that the sequence of functions $(v,\dot v)(\overline x(t_n)+\cdot,
t_n+\cdot)$ converges in the space $\CC^0([0,T],H^1_{\loc}(\Rm)\times
L^2_{\loc}(\Rm))$ to some limit $(w,\dot w)$ which satisfies
\eqref{eq-u}. Note that $|w(0,0)| = \epsilon_0$, because $|v(\overline
x(t_n),t_n)| = \epsilon_0$ for all $n$ by definition of the invasion
point.

Let $c_-, c_+$ be the parabolic speeds associated to $s_-,s_+$
according to \eqref{eq-s-c} (if $s_- \le 0$, we simply take 
$c_- = 0$). We choose any $c \in (c_-,c_+)$ such that $c > c_+-\eta$, 
where $\eta$ is the positive constant given by 
Proposition~\ref{prop-estim-reste}. By construction, $\overline 
y_c(t_n) \to +\infty$ and $\overline y_c(t'_n) \to -\infty$ as 
$n \to \infty$. Applying \eqref{Rcest2} with $y_0 = t_0 = 0$, 
we see that there exists $C_8 > 0$ such that $|R_c(0,t)| \le 
C_8\,e^{-\mu t/2}$ for all $t \ge 0$. Since $\partial_t E_c(0,t)
\le R_c(0,t)$, it follows that $E_c(0,t) \ge E_c(0,t'_n)-  C_9 \,
e^{-\mu t/2}$ for $t \in [0,t_n']$, where $C_9 = 2C_8/\mu$. 
Taking the limit $n \to \infty$ and using \eqref{eq-lowerbound}, 
we conclude that $E_c(0,t) \ge -C_9\,e^{-\mu t/2}$ for all $t \ge 0$. 

On the other hand, using \eqref{def-energy} and the estimate 
above on $R_c(0,t)$, we obtain
\[ 
  E_c(0,t)-E_c(0,0) \,\le\, -(1+\alpha c^2)\int_0^t \int_\Rm e^{cy}
  |\dot v_c(y,\tau)|^2\dd y \dd \tau + C_9 (1-e^{-\mu t/2})~.
\]
Recalling that $E_c(0,t) \ge -C_9 e^{-\mu t/2}$ and setting
$t = t_n + T$, we find
\begin{align*}
  E_c(0,0) &\,\ge\, -C_9 +(1+\alpha c^2)\int_0^{t_n+T} 
  \!\int_\Rm e^{cy}|\dot v_c(y,t)|^2 \dd y \dd t\\
  &\,\ge\, -C_9 +(1+\alpha c^2) \int_{t_n}^{t_n+T} \!\int_\Rm  
  e^{c(\overline y_c(t_n)+y)}|\dot v_c(\overline y_c(t_n)+y,t)|^2 \dd y 
  \dd t~,
\end{align*}
hence
\begin{equation}\label{dem-lemme-asympt-speed}
  \int_0^T \!\int_\Rm  e^{cy}|\dot v_c(\overline y_c(t_n)+y,
  t_n+t)|^2 \dd y \dd t \,\le\, \frac{e^{-c\overline y_c(t_n)}}
  {(1+\alpha c^2)} (E_c(0,0)+C_9) \,\xrightarrow[n\to \infty]{}\, 
  0~.
\end{equation}
Since
\[
  \dot v_c(y,t) \,=\, \dot v\Bigl(\frac {y+ct}{\sqrt{1+\alpha c^2}}\,,
  \,t\Bigr) + sv'\Bigl(\frac {y+ct}{\sqrt{1+\alpha c^2}}\,,
  \,t\Bigr)~, \quad \hbox{where}\quad s = \frac{c}{\sqrt{1+\alpha c^2}}~,
\]
it follows from \eqref{dem-lemme-asympt-speed} that $\int_0^T
\int_{-L}^L |\dot v + sv'|^2(\bar x(t_n)+x,t_n+t)\dd x\dd t$ converges
to zero as $n \to \infty$ for any $L > 0$. Passing to the limit, we
conclude that $\dot w(x,t) + sw'(x,t) = 0$ for all $t \in [0,T]$ and
(almost) all $x \in \Rm$. The key point is that this identity must
hold for all $c \in (c_-,c_+)$ such that $c > c_+-\eta$; i.e., for all
$s$ in a nonempty open interval. Obviously, this implies that $w' =
\dot w = 0$, hence, since $|w(0,0)| = \epsilon_0$, $w(x,t)$ is
identically equal either to $\epsilon_0$ or to $-\epsilon_0$. But this
is impossible, because $w$ must be a solution of \eqref{eq-u} and we 
know from \eqref{beta0deff} that $V'(\pm \epsilon_0) \neq 0$.
\end{demo}


\section{Control of the energy around the invasion point}
\label{sec6}

Proposition~\ref{prop-asympt-speed} shows that the invasion point 
$\overline x(t)$ has an average speed $s_\infty \in 
(0,1/\sqrt{\alpha})$ as $t \to +\infty$. Our next objective 
is to prove that the solution $v(x,t)$ of \eqref{eq-v} 
converges in any neighborhood of the invasion point to the 
profile of a travelling front. To achieve this goal, 
a crucial step is to control the energy $E_c(\overline y_c(t),t)$
for $c$ close to $c_\infty$, where
\begin{equation}\label{cinftydef}  
  c_\infty \,=\, \frac{s_\infty}{\sqrt{1-\alpha s_\infty^2}}~.
\end{equation}
The main result of this section is:

\begin{prop}\label{prop-control-energy}
There exists a positive constant $\eta$ such that, for all
$c\in[c_\infty-\eta,c_\infty+\eta]$, the energy $E_c(\overline
y_c(t),t)$ is a bounded function of $t \ge 0$.
\end{prop}

\begin{rems}\\
{\bf 1.} It is important to realize that, unlike in the previous 
sections, we do not consider in Proposition~\ref{prop-control-energy} 
the energy $E_c(y_0,t)$ located at some fixed point $y_0$ in the 
moving frame, but the energy $E_c(\overline y_c(t),t)$ located 
{\em at the invasion point}. From \eqref{eq-lowerbound} we know 
that $E_c(\overline y_c(t),t) \ge (\min V)/c$ for all $t \ge 0$, 
hence the only problem is to find an upper bound. If $c$ is close 
to $c_\infty$, it follows from Proposition~\ref{prop-estim-reste} 
(with $c_+ = c_\infty$) that $E_c(\overline y_c(t_0),t)$ is bounded 
from above for all $t \ge t_0 \ge 0$. Thus, using the relation
\begin{equation}\label{twopoints}
  E_c(\overline y_c(t),t) \,=\, e^{c(\overline y_c(t_0)-\overline y_c(t))}
  \,E_c(\overline y_c(t_0),t)~, \quad t \ge t_0~,
\end{equation}
which follows immediately from \eqref{eq-Ec-y0}, we see that
$E_c(\overline y_c(t),t)$ is bounded from above if $\overline y_c(t)$
stays bounded from below. Unfortunately $\overline y_c(t) \to -\infty$
as $t \to +\infty$ if $c > c_\infty$, and even if $c = c_\infty$ we do
not know a priori if $\overline y_c(t)$ is bounded from below.  
The essential ingredients in the proof of 
Proposition~\ref{prop-control-energy} are Lemma~\ref{lem-anti-stall},
which allows to control the growth of the exponential factor
$e^{c(\overline y_c(t_0)-\overline y_c(t))}$ in the right-hand side
of \eqref{twopoints}, and Lemma~\ref{lem-dec-energy}, which shows
that the energy $E_c(\overline y_c(t_0),t)$ decays significantly
under appropriate conditions. \\[1mm]
{\bf 2.}  We shall prove in Section~\ref{sec7} that $c_\infty = c_*$
and that the function $v_{c_*}(\overline y_{c_*}(t)+\cdot,t)$
converges uniformly on compact sets to the unique solution $h$ of
\eqref{hdef} such that $h(0) = \epsilon_0$. Now, it is easy to verify
that $h \in H^1_c(\Rm)$ for $c < c_h \equiv \frac12 (c_* +
\sqrt{c_*^2+ 4V''(0)})$. In agreement with
Proposition~\ref{prop-control-energy}, we thus expect that the energy
$E_c(\overline y_c(t),t)$ stays bounded for all times if $c$ is close
to $c_\infty = c_*$, and blows up if $c > c_h$.
\\[1mm]
{\bf 3.}  The fact that the conclusion of
Proposition~\ref{prop-control-energy} holds not only for $c =
c_\infty$ but for all $c$ in a neighborhood of the invasion speed is
one of the key points of our convergence proof. It will allow in
Section~\ref{sec7} to control the variation of the energy
$E_c(\overline y_c(t),t)$ when the parameter $c$ is increased, a
difficult task due to the exponential weight $e^{cy}$ in
\eqref{def-energy}. This problem is completely avoided in the
alternative approach of \cite{Ri} where only bounded weights are used.
\end{rems}

The first step in the proof of Proposition~\ref{prop-control-energy}
consists in showing that the invasion point $\bar x(t)$ cannot make 
arbitrarily large jumps to the left. 

\begin{lemma}\label{lem-anti-stall}
There exists $\eta>0$ such that, for all $c \in (c_\infty-\eta,
c_\infty)$, there exists a positive constant $M_c$ such that,
\begin{equation}\label{eq-anti-stall}
  \overline y_c(t') \,\ge\, \overline y_c(t)-M_c~, \quad
  \hbox{for all } t' \ge t \ge 0~.
\end{equation}
\end{lemma}

\begin{demo}
Let $\eta$ be the positive constant given by Proposition
\ref{prop-estim-reste} for $c_+=c_\infty$. To prove
\eqref{eq-anti-stall}, we argue by contradiction. Assume that there
exist a speed $c \in (c_\infty-\eta,c_\infty)$ and two sequences of
times $\{t_n\}$ and $\{t'_n\}$ such that $t'_n > t_n \ge 0$ for all $n
\in \Nm$ and $\overline y_c(t'_n) - \overline y_c(t_n) \to -\infty$ as
$n \to \infty$. Since $c < c_\infty$, we know that $\overline y_c(t)
\to +\infty$ as $t \to +\infty$, and therefore we must have $t_n \to
+\infty$ as $n \to \infty$. Thus, if we fix any $T > 0$, we 
can apply Proposition~\ref{propcompact} and assume without loss of
generality that the sequence of functions $(v,\dot v)(\overline
x(t_n)+\cdot,t_n+\cdot)$ converges in the space $\CC^0([-T,0],
H^1_\loc(\Rm) \times L^2_\loc(\Rm))$ toward some limit $(w,\dot w)$ 
which satisfies \reff{eq-u}.

Using \eqref{energy-decay} and Proposition~\ref{prop-estim-reste}, 
we obtain for all $n \in \Nm$:
\begin{align*}
  &(1+\alpha c^2)\int_{t_n-T}^{t'_n} \int_\Rm e^{cy} |\dot v_c|^2
  (\overline y_c(t_n)+y,t)\dd y\dd t \\
  &\quad \,\le\, E_c(\overline y_c(t_n),t_n-T) 
  - E_c(\overline y_c(t_n),t'_n) + \frac{2K_4}{\mu} \max\Bigl(
  E_c(\overline y_c(t_n),t_n-T)\,,\, e^{-c\overline y_c(t_n)}\Bigr)~.
\end{align*} 
In view of \eqref{twopoints}, we have
\[
  \max\Bigl(E_c(\overline y_c(t_n),t_n-T)\,,\,
  e^{-c\overline y_c(t_n)}\Bigr) \,=\, e^{-c\overline y_c(t_n)}
  \max(E_c(0,t_n-T)\,,1) \,\xrightarrow[n\to \infty]{}\, 0~,
\]
because $\overline y_c(t_n) \to +\infty$ and $E_c(0,t_n-T)$ 
is bounded from above due to \eqref{Ecbd1}. On the other hand, 
since $\overline y_c(t'_n) - \overline y_c(t_n) \to -\infty$
by assumption, it follows from \eqref{eq-lowerbound} that
\[
  - E_c(\overline y_c(t_n),t'_n) \,\le\, e^{c(\overline y_c(t'_n)
  -\overline y_c(t_n))}\,\frac{|\min V|}{c} \,\xrightarrow[n\to 
  \infty]{}\, 0~.
\]
Thus, we have shown:
\[
  (1+\alpha c^2)\int_{-T}^0\int_\Rm e^{cy} |\dot v_c|^2(\overline 
  y_c(t_n)+y,t_n+t)\dd y\dd t \,\xrightarrow[n\to \infty]{}\, 0~.
\]

Proceeding as in the proof of Lemma~\ref{lemme-asympt-speed}, we
conclude that $\dot w(x,t) + sw'(x,t) = 0$ for all $t \in [-T,0]$ and
(almost) all $x \in \Rm$, where $s = c/\sqrt{1+\alpha c^2}$.  Now, the
crucial observation is that, for any $c' \in (c,c_\infty)$, we still
have $\overline y_{c'}(t_n) \to +\infty$ and $\overline y_{c'}(t_n') -
\overline y_{c'}(t_n) \to -\infty$ as $n \to \infty$. The second
claim follows immediately from the identity
\begin{equation}\label{baryid}
  \overline y_{c_2}(t') - \overline y_{c_2}(t) \,=\, 
  \frac{\sqrt{1+\alpha c_2^2}}{\sqrt{1+\alpha c_1^2}}
  \,\Bigl(\overline y_{c_1}(t') - \overline y_{c_1}(t)\Bigr)
  + \Bigl(\frac{\sqrt{1+\alpha c_2^2}}{\sqrt{1+\alpha c_1^2}}
  \,c_1 - c_2\Bigr)(t'-t)~. 
\end{equation}
Thus, repeating the same arguments, we conclude that $\dot w + 
s'w' = 0$ for all $s'$ in a nonempty open interval. This implies that
$\dot w = w' = 0$, and we obtain a contradiction as in 
Lemma~\ref{lemme-asympt-speed}. 
\end{demo}

Combining the bound \eqref{eq-anti-stall} and the identity 
\eqref{baryid}, we obtain the following useful estimate, 
which is valid in any frame whose speed is close enough 
to the invasion speed. 

\begin{coro}\label{barycor} 
For all $p>0$, there exist $\eta > 0$ and $M > 0$ such that, for all
$c \in [c_\infty-\eta, c_\infty+\eta]$,
\begin{equation}\label{lowbdd}
  \overline y_c(t') - \overline y_c(t) \,\ge\, -M -p(t'-t)~, 
  \quad \hbox{for all } t' \ge t \ge 0~.
\end{equation}
\end{coro}

The next proposition shows that, if the energy is sufficiently large at
a given time, and if the invasion point stays bounded from above on a
sufficiently long time interval, then a significant decay of energy
must occur. 

\begin{lemma}\label{lem-dec-energy}
There exist positive constants $\eta$, $t_0$, $\omega$, and 
$K_5$ such that the following holds. For any $c \in [c_\infty-\eta,
c_\infty+\eta]$, any $y_0 \in \Rm$, any $t_1 \ge t_0$, any 
$T > 0$, and any $M > 0$, if the invasion point satisfies 
$\overline y_c(t) \le y_0 + M$ for all $t \in [t_1,t_1+T]$, then 
\begin{equation}\label{newdec-energy}
  E_c(y_0,t_1+T) \,\le\, K_5 \Bigl(e^{-\omega T}\,E_c(y_0,t_1)
  + e^{cM}\Bigr)~.
\end{equation}
\end{lemma}

\begin{demo}
Let $\eta$ be the positive number given by Proposition
\ref{prop-estim-reste} for $c_+=c_\infty$ and let
$c\in[c_\infty-\eta,c_\infty+\eta]$. Given $y_0 \in \Rm$, we 
define
\[
  \EE_c(y_0,t) \,=\, \int_\Rm e^{cy}\Bigl(\frac\alpha 2|\dot
  v_c|^2 + \frac12 |v'_c|^2 + V(v_c) + \alpha \gamma v_c\dot
  v_c\Bigr)(y_0+y,t)\dd y~,
\]
where $\gamma > 0$ will be fixed later. Equation \eqref{eq-vc-rc} 
satisfied by $v_c$ implies that
\begin{align*}
  \partial_t \EE_c(y_0,t) =\int_{\Rm} e^{cy}\Bigl(&-(1+\alpha c^2)
  |\dot v_c|^2 + f(v_c,r_c)\dot v_c + \alpha \gamma |\dot v_c|^2 - 
  \gamma(1+2\alpha c^2)v_c \dot v_c \\
  &- 2\alpha\gamma c v_c'\dot v_c - \gamma|v'_c|^2- \gamma 
  V'(v_c)v_c + \gamma f(v_c,r_c)v_c \Bigr)(y_0+y,t)\dd y~.
\end{align*}
From \eqref{fbound} and \eqref{vcrcbound} we know that 
\[
  \int_{\Rm} e^{cy} |f(v_c,r_c)| (|\dot v_c| + \gamma |v_c|)
  (y_0+y,t)\dd y \,\le\, K_2 \epsilon_0 \,e^{-\mu t} \int_{\Rm} 
  e^{cy}(|v_c|^2 + |\dot v_c|^2)(y_0+y,t)\dd y~,
\]
provided that $\gamma \le 3/4$. Thus, using the bound $2ab \le 
C^{-1} a^2 + C b^2$ and the Poincar\'e inequality \eqref{eq-Poincare}, 
we obtain that, if $\gamma$ is small enough and $t_0$ is large enough, 
the following estimate holds for all $t \ge t_0$:
\[
  \partial_t \EE_c(y_0,t) \,\le\, -\int_{\Rm} e^{cy}\Bigl(
  \frac12 |\dot v_c|^2 + \frac\gamma2 |v'_c|^2 + \gamma V'(v_c)
  v_c\Bigr)(y_0+y,t)\dd y~.
\]
On the other hand, we know from \eqref{vcrcbound} that $|v_c(y,t)|$ 
is uniformly bounded for all $y \in \Rm$ and all $t \ge 0$, and 
from \eqref{beta0deff} that $2V'(v_c(y,t))v_c(y,t) \ge V(v_c(y,t)) 
\ge (\beta_0/4)v_c(y,t)^2$ for all $y \ge \overline y_c(t)$. Thus, 
there exist $\omega > 0$ and $C_0 > 0$ such that
\begin{eqnarray}\nonumber
  \partial_t \EE_c(y_0,t) &\le& -\omega \EE_c(y_0,t) + C
  \int_{-\infty}^{\overline y_c(t)-y_0} e^{cy} \Bigl(|V'(v_c)v_c| + 
  |V(v_c)|\Bigr)(y_0+y,t)\dd y~,\\ \label{dem-prop-dec-en}
  &\le& -\omega \EE_c(y_0,t) + C_0 \,e^{c(\overline y_c(t)-y_0)}~,
\end{eqnarray}
for all $t \ge t_0$. In a similar way, there exist $C_1 > 1$ 
and $C_2 > 0$ such that 
\begin{equation}\label{equivenerg}
  C_1^{-1}E_c(y_0,t) - C_2 \,e^{c(\overline y_c(t)-y_0)} \,\le\, 
  \EE_c(y_0,t) \,\le\, C_1 E_c(y_0,t) + C_2 \,e^{c(\overline y_c(t)-y_0)}~,
\end{equation}
for all $t \ge t_0$. Remark that, in \eqref{dem-prop-dec-en} and 
\eqref{equivenerg}, all constants can be chosen to be independent 
of $y_0 \in \Rm$, of $c \in [c_\infty-\eta,c_\infty+\eta]$, and 
of $t \ge t_0$.

Now, we fix $t_1 \ge t_0$ and assume that $\overline y_c(t)-y_0 \le 
M$ when $t \in [t_1,t_1+T]$, for some $M > 0$ and some $T > 0$. 
Integrating the differential inequality \eqref{dem-prop-dec-en}, 
we find
\[
  \EE_c(y_0,t_1+T) \,\le\, e^{-\omega T}\EE_c(y_0,t_1) + 
  \frac{C_0}{\omega}\,e^{cM}~.
\]
Combining this result with \eqref{equivenerg}, we arrive at
\[
  E_c(y_0,t_1+T) \,\le\, C_1^2 \,e^{-\omega T}E_c(y_0,t_1)
  + \Bigl(2C_1 C_2 + \frac{C_0 C_1}{\omega}\Bigr)\,e^{cM}~,
\]
which is the desired estimate. 
\end{demo}

Using the control on the invasion point given by 
Lemma~\ref{lem-anti-stall} and the decay of energy 
described in Lemma~\ref{lem-dec-energy}, we are now able 
to prove the main result. 

\medskip\noindent 
{\textbf{Proof of Proposition~\ref{prop-control-energy}:}} 
Let $t_0$, $\omega$, $K_5$ be as in Lemma~\ref{lem-dec-energy}, 
and choose $p > 0$ such that $4 c_\infty p \le \min(\omega,\mu)$, 
where $\mu > 0$ is as in Lemma~\ref{lemme-estim-reste}. 
By Corollary~\ref{barycor}, there exist $\eta > 0$ and $M > 0$ 
such that \eqref{lowbdd} holds, and without loss of generality
we can assume that $\eta < c_\infty$ and that $M$ is large enough 
so that $e^{-(c_\infty-\eta)M} \le 1/2$. In the rest of the proof, 
we fix some $c \in [c_\infty-\eta,c_\infty+\eta] \subset 
(0,2c_\infty)$ (but all constants will be independent of $c$).

From \eqref{energy-decay} and Proposition~\ref{prop-estim-reste}
we know that, for all $t \ge 0$ and all $\tau \ge 0$,  
\begin{align*}
  \partial_\tau E_c(\overline y_c(t),t+\tau) &\,\le\, K_4 
  \,e^{-\mu (t+\tau)/2} \,\max\Bigl(E_c(\overline y_c(t),t)\,,\,
    e^{-c\overline y_c(t)}\Bigr)\\
  &\,\leq\, K_4 \,e^{-\mu (t+\tau)/2} \,e^{-c\overline y_c(t)}
  \,\max\Bigl(E_c(0,t)\,,\,1\Bigr)~.
\end{align*}
Since, by \eqref{lowbdd}, $c\overline y_c(t) \ge c(\overline
y_c(0)-M-pt) \ge -C - \mu t/2$ for all $t \ge 0$, and since $E_c(0,t)$
is uniformly bounded from above by Proposition~\ref{prop-estim-reste},
there exists $C_3 > 0$ such that, for all $t \ge 0$ and all $\tau 
\ge 0$,
\begin{equation}\label{dem-th-control-energy}
  E_c(\overline y_c(t),t+\tau) \,\le\, E_c(\overline y_c(t),t) + C_3~.
\end{equation}
Now, we choose $T \ge 1$ and $C_4 \ge C_3$ such that the following
inequalities hold: 
\[
  4 K_5 \,e^{cM} \,e^{-\omega T/2} \,\le\, 1~, \quad 
  \hbox{and} \quad 4 K_5 \,e^{4cM} \,e^{cp(1+T)} \,\le\, C_4~.
\]
We claim that, if $E_c(\overline y_c(t),t) \ge C_4$ for some 
$t \ge t_0$, then there exists $t' \in [1,T]$ such that 
$E_c(\overline y_c(t+t'),t+t') \le \frac12 E_c(\overline y_c(t),t)$.

To prove this claim, we distinguish two possible cases. If 
there exists $\tau \in [0,T]$ such that $\overline y_c(t+\tau) \ge
\overline y_c(t) + 3M + p$, then $t'=\max (\tau,1)$ is a suitable 
choice. Indeed, by Corollary~\ref{barycor}, we have $\overline 
y_c(t+t')\ge \overline y_c(t) + 2M$ and so, using 
\eqref{dem-th-control-energy}, we find
\begin{align*}
  E_c(\overline y_c(t+t'),t+t') &\,=\, e^{c(\overline y_c(t)-\overline
    y_c(t+t'))} E_c(\overline y_c(t),t+t')\\
  &\,\le\, e^{-2cM} (E_c(\overline y_c(t),t) + C_3)
  \,\le\, {\textstyle \frac12} E_c(\overline y_c(t),t)~,
\end{align*}
because $e^{-2cM} \le 1/4$ and $C_3\leq C_4\le E_c(\overline y_c(t),t)$ 
by assumption. On the other hand, if $\overline y_c(t+\tau) \le  
\overline y_c(t) + 3M + p$ for all $\tau \in [0,T]$, we can 
take $t' = T$ because, due to Corollary~\ref{barycor},
Lemma~\ref{lem-dec-energy} and our choices of $T$, $C_4$, and $p$, 
we have
\begin{align*}
  E_c(\overline y_c(t+T),t+T) &\,=\, e^{c(\overline y_c(t)-\overline
    y_c(t+T))} E_c(\overline y_c(t),t+T)\\
  &\,\le\, e^{c(M+pT)} K_5 \Bigl(e^{-\omega T} E_c(\overline y_c(t),t)
  + e^{c(3M+p)}\Bigr) \,\le\, {\textstyle \frac12} 
  E_c(\overline y_c(t),t)~.
\end{align*}

We now show that the claim above implies 
Proposition~\ref{prop-control-energy}. To this purpose, we construct 
the following sequence of times. We take $t_0 > 0$ as in 
Lemma~\ref{lem-dec-energy}, and given $t_n$ we define 
$t_{n+1}$ in the following way. If $E_c(\overline y_c(t_n),t_n)\le
C_4$, we simply set $t_{n+1} = t_n+1$ and, using 
\eqref{dem-th-control-energy} and Corollary~\ref{barycor}, we get 
\[
  E_c(\overline y_c(t_{n+1}),t_{n+1}) \,=\, e^{c(\overline y_c(t_{n})
  -\overline y_c(t_{n+1}))} E_c(\overline y_c(t_{n}),t_{n+1}) 
  \,\le\, e^{c(M+p)}(C_4 + C_3)~.
\]
If $E_c(\overline y_c(t_n),t_n)\ge C_4$ we set $t_{n+1} = t_n+t'$,
where $t' \ge 1$ is the time given by the claim above when $t=t_n$.
In this case, we know that $E_c(\overline y_c(t_{n+1}),t_{n+1})\le
\frac12 E_c(\overline y_c(t_n),t_n)$. By contruction, the sequence
$\{t_n\}$ goes to $+\infty$ as $n \to \infty$, and 
\[
  \limsup_{n\to \infty} E_c(\overline y_c(t_n),t_n) \,\le\, 
  e^{c(M+p)}(C_4 + C_3)~.
\]
To conclude the proof, we observe that \eqref{dem-th-control-energy}
and Corollary~\ref{barycor} provides a control on the energy for the
remaining times. Indeed, since $t_{n+1}-t_n \le T$ for all $n$, we
have $\overline y_c(t)\ge
\overline y_c(t_n)-M-pT$ for all $t\in[t_n,t_{n+1}]$, hence  
\[
  E_c(\overline y_c(t),t) \,\le\, e^{c(M+pT)}E_c(\overline y_c(t_n),t)
  \,\le\, e^{c(M+pT)}(E_c(\overline y_c(t_n),t_n) +C_3)~, \quad
  t \in [t_n,t_{n+1}]~.
\]
A similar argument shows that $E_c(\overline y_c(t),t)$ is bounded
from above for $t \in [0,t_0]$. \hfill $\square$


\section{Convergence to a travelling wave}
\label{sec7}

The purpose of this section is to show that, for any $L > 0$, the
solution $v$ of \eqref{eq-v} converges to a travelling front uniformly
in the interval $(\overline x(t)-L,+\infty)$. The key step is to prove
that, in the frame moving at the invasion speed $s_\infty$, the energy
dissipation around the invasion point converges to zero as $t \to
+\infty$.

\begin{prop}\label{prop-relaxation}
Let $s_\infty$ be the invasion speed introduced in 
Proposition~\ref{prop-asympt-speed} and let $c_\infty$ be the 
parabolic speed \eqref{cinftydef}. For any $T > 0$, we 
have
\[
  \int_{t-T}^{t}\int_\Rm e^{c_\infty y}|\dot v_{c_\infty}(\overline
  y_{c_\infty}(t)+y,\tau)|^2 \dd y \dd\tau 
  \,\xrightarrow[t\to +\infty]{}\, 0~.
\] 
\end{prop}

We start the proof with an auxiliary result showing that the 
energy $E_c(\overline y_c(t),t)$ is a continuous function of 
the parameter $c$.

\begin{lemma}\label{lipschitz}
Let $\eta > 0$ be as in Proposition~\ref{prop-control-energy}. Given any
$T \ge 0$, there exists $K_6 > 0$ such that, for all $c_1,c_2 
\in [c_\infty-\eta/2,c_\infty+\eta/2]$, all $t \ge T$ and all 
$\tau \in [t-T,t]$, the following estimate holds:
\[
  |E_{c_1}(\overline y_{c_1}(t),\tau) - 
   E_{c_2}(\overline y_{c_2}(t),\tau)| \,\le\, K_6 |c_1-c_2|~.
\]
\end{lemma}

\begin{demo}
If we return to the original variables using the definitions 
\eqref{def-vc-rc}, \eqref{barydef}, we obtain the identity
\begin{align}\label{iden}
  &E_c(\overline y_c(t),\tau) \,=\, \int_\Rm e^{cy}\Bigl(\frac\alpha2 
  |\dot v_c|^2 + \frac12 |v'_c|^2 + V(v_c) \Bigr)(\overline y_c(t)
  +y,\tau)\dd y \\ \nonumber
  &\quad \,=\, \sqrt{1{+}\alpha c^2}\,e^{c(\overline y_c(\tau)-
  \overline y_c(t))} \int_\Rm e^{c\sqrt{1+\alpha c^2}x} \Bigl(\frac\alpha2 
  |\dot v {+} sv'|^2 + \frac{|v'|^2}{2(1{+}\alpha c^2)} + V(v)\Bigr)
  (\overline x(\tau){+} x,\tau)\dd x\,,
\end{align}
where $s = c/\sqrt{1+\alpha c^2}$. Assume first that $\tau = t$ and 
$c = \overline c$, where $\overline c = c_\infty + \eta$. We know 
from Proposition~\ref{prop-control-energy} that $E_{\overline c}(\overline 
y_{\overline c}(t),t)$ is bounded (from above) for all times.  
On the other hand, we obtain a lower bound on the last member of 
\eqref{iden} if we replace $V(v(\overline x(t)+x,t))$ by zero 
if $x \ge 0$ and by $\min V$ if $x \le 0$. Thus, using in 
addition the Poincar\'e inequality \eqref{eq-Poincare} to control
$|v|^2$ in terms of $|v'|^2$, we deduce from \eqref{iden} 
that there exists a constant $C_0 > 0$ such that
\begin{equation}\label{unifbound}
  \int_\Rm e^{\overline c \sqrt{1+\alpha \overline c^2}x} \Bigl(|\dot v|^2
  +|v'|^2 + |v|^2\Bigr)(\overline x(t)+x,t)\dd x \,\le\, C_0~, 
  \quad \hbox{for all } t \ge 0~.
\end{equation}
Using the uniform control
\eqref{unifbound}, it is a straightforward exercise to verify that the
last member of \eqref{iden} is indeed a Lipschitz function of $c \in
[c_\infty-\eta/2,c_\infty+\eta/2]$, uniformly in $t \ge T$ and $\tau
\in [t-T,t]$. The only potential difficulty comes from the exponential 
terms. If we denote
\[
  Y_c(t,\tau) \,=\, c(\overline y_c(\tau)-\overline y_c(t)) 
  \,=\, c\sqrt{1+\alpha c^2}\,(\overline x(\tau)-\overline x(t))
  + c^2(t-\tau)~, 
\]
we know from Corollary~\ref{barycor} that $Y_c(t,\tau) \le 
c(M+pT)$ for all $\tau \in [t-T,t]$ and all $c \in [c_\infty-\eta,
c_\infty+\eta]$, hence
\[
  \Bigl|e^{Y_{c_1}(t,\tau)} - e^{Y_{c_2}(t,\tau)}\Bigr| \,\le\, 
  e^{\max(Y_{c_1}(t,\tau),Y_{c_2}(t,\tau))} \,|Y_{c_1}(t,\tau) - 
  Y_{c_2}(t,\tau)| \,\le\, C_1 |c_1 - c_2|~.
\]
On the other hand, if $c_1,c_2 \in [c_\infty-\eta/2,c_\infty+\eta/2]$,
we can bound
\[
  \Bigl|e^{c_1\sqrt{1+\alpha c_1^2}x} - e^{c_2\sqrt{1+\alpha c_2^2}x}\Bigr|
  \,\le\, C_2 |c_1 - c_2| \left\{
  \begin{array}{rcl} 
  e^{\overline c\sqrt{1+\alpha \overline c^2}x} & \hbox{if} & x \ge 0~,\\
  e^{\underline c\sqrt{1+\alpha \underline c^2}x} & \hbox{if} & x \le 0~,
  \end{array}\right.
\]
where $\overline c = c_\infty+\eta$ and $\underline c = c_\infty-\eta$.
Thus, using estimate \eqref{unifbound} for $x \ge 0$ and the
uniform bound \eqref{vbound} for $x \le 0$, we obtain the desired 
conclusion.
\end{demo}

\medskip\noindent 
{\textbf{Proof of Proposition~\ref{prop-relaxation}:}}
We argue by contradiction. Assume that there exist $\hat\delta > 0$, 
$T > 0$, and a sequence of times $\{t_n\}$ going to $+\infty$ 
such that
\begin{equation}\label{contradict}
  \int_{t_n-T}^{t_n}\int_\Rm e^{c_\infty y} \,|\dot v_{c_\infty}
  (\overline y_{c_\infty}(t_n)+y,t)|^2 \dd y\dd t \,\ge\, \hat\delta~,
\end{equation}
for all $n \in \Nm$. Following an idea introduced in 
\cite{Ri}, we shall arrive at a contradiction by considering 
the variation of the energy along the {\em broken line} 
connecting the points $(t_n,\overline y_{c_\infty}(t_n))$ in 
the $(t,y)$ plane, see Fig.~2.

\begin{figure}[hbt]
\vspace{0.4cm}
\begin{center}
\hspace{-1cm}\input{fig2.pstex_t} 
\end{center}
\vskip 0.0truecm{\baselineskip=0.8\baselineskip \noindent \vbox{\noindent
{\footnotesize {\bf Fig.~2:} The behavior of the energy along the
broken line.}}\par}
\vspace{0.4cm}
\end{figure}

To this end we define, for all $n \in \Nm$,
\[
  s_n \,=\,\frac 1{\sqrt{1+\alpha c_\infty^2}}\frac {\overline
  y_{c_\infty}(t_{n+1})- \overline y_{c_\infty}(t_{n})}
  {t_{n+1}-t_n} + s_\infty~, \quad \hbox{and}\quad 
  c_n \,=\, \frac{s_n}{\sqrt{1-\alpha s_n^2}}~.
\]
As is easily verified, the speed $s_n$ is the slope of the line
segment connecting $(t_n,\overline x(t_n))$ and $(t_{n+1},\overline
x(t_{n+1}))$; i.e., $\overline x(t_{n+1}) - \overline x(t_n) =
s_n(t_{n+1}-t_n)$. Since $\overline y_{c_\infty}(t)/t$ converges to
zero as $t \to +\infty$ by Proposition~\ref{prop-asympt-speed}, we 
can assume (up to extracting a subsequence) that $s_n \in 
(0,1/\sqrt{\alpha})$ for all $n \in \Nm$ and that $s_n \to s_\infty$
as $n \to \infty$. Then the parabolic speed $c_n$ is well-defined
for all $n$, and by construction $\overline y_{c_n}(t_{n+1}) = 
\overline y_{c_n}(t_{n})$. Extracting another subsequence if 
needed, we can further assume that $t_{n+1}\ge t_n+T$ and 
$|c_n - c_\infty| \le \eta/2$ for all $n \in \Nm$, where 
$\eta > 0$ is as in Proposition~\ref{prop-control-energy}, and 
that the sum $\sum_{n\ge0} |c_n - c_\infty|$ is finite. 

Now we define, for each $n \in \Nm$,
\begin{eqnarray*}
  \Delta_n &=& E_{c_n}(\overline y_{c_n}(t_n),t_n) - 
  E_{c_{n+1}}(\overline y_{c_{n+1}}(t_{n+1}),t_{n+1}) \\
  &=& E_{c_n}(\overline y_{c_n}(t_n),t_n) - 
  E_{c_n}(\overline y_{c_n}(t_{n+1}),t_{n+1}) \\
  &+& E_{c_n}(\overline y_{c_n}(t_{n+1}),t_{n+1}) - 
  E_{c_{n+1}}(\overline y_{c_{n+1}}(t_{n+1}),t_{n+1}) 
  \,=\, \Delta_n^1 + \Delta_n^2~.
\end{eqnarray*}
Since $\overline y_{c_n}(t_{n+1}) = \overline y_{c_n}(t_{n})$, 
the quantity $\Delta_n^1$ is the variation of the energy 
$E_{c_n}(y_0,t)$ at a fixed point $y_0 \in \Rm$ on the time
interval $[t_n,t_{n+1}]$. By \eqref{energy-decay} and 
Proposition~\ref{prop-estim-reste}, we have
\[
  \Delta_n^1 \,=\, (1+\alpha c_n^2)\int_{t_n}^{t_{n+1}}\!\!\int_\Rm
  e^{c_n y}|\dot v_{c_n}(\overline y_{c_n}(t_{n+1})+y,t)|^2 \dd y
  \dd t - \int_{t_n}^{t_{n+1}} R_{c_n}(\overline y_{c_n}(t_n),t)
  \dd t~,
\]
and
\[
  |R_{c_n}(\overline y_{c_n}(t_n),t)| \,\le\, K_4 \,e^{-\mu t/2}
  \,e^{-c_n \overline y_{c_n}(t_n)} \,\max(E_{c_n}(0,t_n)\,,\,1)~.
\]
But $E_{c_n}(0,t_n)$ is bounded by from above uniformly in $n$ by
Proposition~\ref{prop-estim-reste}, and since $\overline y_{c_n}
(t_n)/t_n$ converges to zero as $n \to \infty$ we can assume 
without loss of generality that $c_n \overline y_{c_n}(t_n) \ge 
-\mu t_n/4$ for all $n \in \Nm$. Thus
\[
  \Delta_n^1 \,\ge\, \int_{t_{n+1}-T}^{t_{n+1}}\int_\Rm
  e^{c_n y}|\dot v_{c_n}(\overline y_{c_n}(t_{n+1})+y,t)|^2 \dd y
  \dd t - C_3 \,e^{-\mu t_n/4}~,
\]
for some $C_3 > 0$. Moreover, since $|c_n - c_\infty| \le 
\eta/2$, the proof of Lemma~\ref{lipschitz} shows that, for 
all $t \in [t_{n+1}-T,t_{n+1}]$,
\[
  \left|\int_\Rm e^{c_n y}|\dot v_{c_n}(\overline y_{c_n}(t_{n+1})
  {+}y,t)|^2 \dd y \dd t - \int_\Rm e^{c_\infty y}|\dot 
  v_{c_\infty}(\overline y_{c_\infty}(t_{n+1}){+}y,t)|^2 
  \dd y \dd t \right| \,\le\, C_4 |c_n-c_\infty|\,,
\]
for some $C_4 > 0$. Combining both estimates and using 
the assumption \eqref{contradict}, we thus obtain
\[
  \Delta_n^1 \,\ge\, \hat\delta - C_4 T |c_n-c_\infty| - 
  C_3 \,e^{-\mu t_n/4}~, \quad n \in \Nm~.
\]
On the other hand, the quantity $\Delta_n^2$ represents the 
change in the energy $E_c(\overline y_c(t_{n+1}),t_{n+1})$ 
when $c$ varies from $c_n$ to $c_{n+1}$. By Lemma~\ref{lipschitz}, 
we have $|\Delta_n^2| \le K_6 |c_n-c_{n+1}|$, hence 
\[
  \Delta_n \,=\, \Delta_n^1 + \Delta_n^2 \,\ge\, \hat\delta 
  - K_6 |c_n-c_{n+1}| - C_4 T |c_n-c_\infty| - C_3 \,e^{-\mu t_n/4}~, 
  \quad n \in \Nm~.
\]

To conclude the proof, we observe that
\[
  E_{c_0}(\overline y_{c_0}(t_0),t_0) - 
  E_{c_N}(\overline y_{c_N}(t_N),t_N) \,=\, \sum_{n=0}^{N-1} 
  \Delta_n \,\xrightarrow[N\to \infty]{}\, +\infty~,
\]
because $t_n \ge nT$ and the sum $\sum_{n\ge0} |c_n - c_\infty|$ is
finite. Thus $E_{c_N}(\overline y_{c_N}(t_N),t_N) \to -\infty$ as $N
\to \infty$, in contradiction with the lower bound
\eqref{eq-lowerbound}. Thus \eqref{contradict} cannot hold for all $n
\in \Nm$, and Proposition~\ref{prop-relaxation} is proved. \hfill$\square$

\medskip
Using now, for the first time, the fact that the differential 
equation \eqref{hdef} has a front-like solution for a single 
value of $c_*$ and that this solution is unique up to translations, 
we can establish the local convergence to a travelling front.

\begin{coro}\label{loc-conv}
The invasion speed satisfies $s_\infty = s_* \equiv c_*/\sqrt{
1+\alpha c_*^2}$, and we have
\begin{align*}
  \int_\Rm e^{c_* \sqrt{1+\alpha c_*^2}\,x} \Bigl(
  |\dot v(\overline x(t)+x,t) + s_* v_*'(x)|^2 &+ 
  |v'(\overline x(t)+x,t) - v_*'(x)|^2 \\ &+ 
  |v(\overline x(t)+x,t) - v_*(x)|^2\Bigr)\dd x 
  \,\xrightarrow[t\to +\infty]{}\, 0~,
\end{align*}
where $v_*(x) = h(\sqrt{1+\alpha c_*^2}\,x)$ and $h$ is the solution
of \eqref{hdef} normalized so that $h(0) = \epsilon_0$. In particular,
$v$ converges to a front uniformly in any interval of the type
$(\overline x(t)-L,+\infty)$.
\end{coro}

\begin{demo}
Fix $T > 0$ and let $\{t_n\}$ be a sequence of times going to 
$+\infty$ as $n \to \infty$. In view of 
Proposition~\ref{propcompact}, we can assume that the sequence of 
functions $(v,\dot v)(\overline x(t_n)+\cdot,t_n+\cdot)$ converges 
in the space $\CC^0([-T,0],H^1_{\loc}(\Rm)\times
L^2_{\loc}(\Rm))$ to some limit $(w,\dot w)$ which satisfies
\eqref{eq-u}. By \eqref{def-vc-rc} and Proposition~\ref{prop-relaxation},
we have 
\[
  \int_{-T}^0\int_\Rm e^{c_\infty \sqrt{1+\alpha c_\infty^2}\,x}
  |\dot v + s_\infty v'|^2(\overline x(t_n)+x,t_n+t)\dd x\dd t
  \,\xrightarrow[n\to \infty]{}\, 0~,
\]
hence $\dot w(x,t) + s_\infty w'(x,t) = 0$ for all $t \in [-T,0]$ and
(almost) all $x \in \Rm$. Setting $w(x,t) = h(\sqrt{1+\alpha
  c_\infty^2}\,x - c_\infty t)$, we see that $h$ is a solution of the
differential equation $h'' + c_\infty h' - V'(h) = 0$. We also know
that $|h(x)| \le M_0$ for all $x \le 0$, that $|h(x)| \le \epsilon_0$
for all $x \ge 0$, and that $|h(0)| = |w(0,0)| = \epsilon_0$. Due to
our assumptions \eqref{coercive}--\eqref{critval} on the potential
$V$, these properties together imply that $c_\infty = c_*$, hence
$s_\infty = s_*$, and that $h$ is the unique solution of \eqref{hdef}
such that $h(0) = \epsilon_0$, see e.g. \cite{AW}. Since the limit is
unique, we observe that the convergence above holds in fact for any
sequence $t_n \to +\infty$. In particular, if we denote $v_*(x) =
h(\sqrt{1+\alpha c_*^2}\,x)$, we conclude that $(v,\dot v) (\overline
x(t)+\cdot,t)$ converges as $t \to +\infty$ to $(v_*,-s_* v_*')$ in
$H^1([-L,L])\times L^2([-L,L])$, for any $L > 0$. Using in addition
the estimate \eqref{unifbound} for $x \ge L$, and the uniform bound
\eqref{vbound} for $x \le -L$, we obtain the desired conclusion.
\end{demo}

One can extract from the proof of Corollary~\ref{loc-conv}
the following useful information on the invasion point:

\begin{lemma}\label{varxbar}
For any $T > 0$ we have
\[
  \sup_{|\tau|\le T} |\overline x(t+\tau) - \overline x(t)
  - s_* \tau| \,\xrightarrow[t\to +\infty]{}\, 0~.
\]
\end{lemma}

\begin{demo}
Fix $T > 0$, and choose $L > 0$ large enough so that
\begin{equation}\label{hinterval}
  h(\sqrt{1+\alpha c_*^2}\,L) \,\le\, \frac{\epsilon_0}2~, \quad 
  \hbox{and} \quad h(-\sqrt{1+\alpha c_*^2}\,L + c_*T) \,\ge\, 
  \frac{1+\epsilon_0}2~,
\end{equation}
where $h$ is as in Corollary~\ref{loc-conv}. We claim that, for 
any $\hat\delta > 0$, there exists $t_0 \ge T$ such that, for all 
$t \ge t_0$, 
\begin{equation}\label{deltaconv}
  \sup_{\tau \in [-T,0]} \,\sup_{x \ge -L} |v(\overline x(t)+x,t+\tau)
  - h(\sqrt{1+\alpha c_*^2}\,x - c_*\tau)| \,\le\, \hat\delta~.
\end{equation}
Indeed, if we restrict the values of $x$ to a bounded interval 
$I = [-L,L']$, where $L' > 0$, the analog of \eqref{deltaconv}
follows immediately from the proof of Corollary~\ref{loc-conv}
and the fact that $H^1(I) \hookrightarrow L^\infty(I)$. On 
the other hand, by Lemma~\ref{lem-anti-stall}, there 
exists $C > 0$ such that $\overline x(t) \ge \overline x(t+\tau) 
-C$ for all $\tau \in [-T,0]$. In view of \eqref{unifbound} 
we thus have
\begin{align*}
  \sup_{\tau \in [-T,0]} &\,\sup_{x \ge L'} |v(\overline x(t)+x,t+\tau) - 
  h(\sqrt{1+\alpha c_*^2}\,x - c_*\tau)| \\
  &\,\le\, \sup_{\tau \in [-T,0]}\,\sup_{x \ge L'-C} |v(\overline 
  x(t+\tau)+x,t+\tau)| + h(\sqrt{1+\alpha c_*^2}\,L')
  \,\xrightarrow[L'\to +\infty]{}\, 0~,
\end{align*}
uniformly in $t$. This proves \eqref{deltaconv}. 

We now assume that $\hat\delta < \min(\epsilon_0,1-\epsilon_0)/2$. For 
any $t \ge t_0$ and any $\tau \in [-T,0]$, it follows from 
\eqref{hinterval}, \eqref{deltaconv} that
\[
  |v(\overline x(t)-L,t+\tau)| \,>\, \epsilon_0~, \quad 
  \hbox{and} \quad \sup_{x \ge L} |v(\overline x(t)+x,t+\tau)| 
  \,<\, \epsilon_0~.
\]
By the definition \eqref{barxdef} of the invasion point, this means
that $\overline x(t+\tau) \in [\overline x(t)-L,\overline x(t)+L]$. 
Using \eqref{deltaconv} with $x = \overline x(t+\tau) - \overline x(t)$ 
and recalling that $v(\overline x(t+\tau),t+\tau) = h(0) =
\epsilon_0$, we obtain
\[
  \hat\delta \,\ge\, |\epsilon_0 - h(\sqrt{1+\alpha c_*^2}\,(\overline 
  x(t+\tau) - \overline x(t)) - c_*\tau)| \,\ge\, m \sqrt{1+\alpha
  c_*^2}\,|\overline x(t+\tau) - \overline x(t) - s_* \tau|~,
\] 
where
\[
  m \,=\, \min\Bigl\{|h'(y)| \,\Big|\, -\sqrt{1+\alpha c_*^2}\,L 
  \le y \le \sqrt{1+\alpha c_*^2}\,L + c_*T\Bigr\} \,>\, 0~.
\]
Thus $|\overline x(t+\tau) - \overline x(t) - s_* \tau| \le 
\hat\delta/(m \sqrt{1+\alpha c_*^2})$ for all $t \ge t_0$ and 
all $\tau \in [-T,0]$. Since $\hat\delta > 0$ was arbitrarily small, 
we obtain the desired conclusion. 
\end{demo}


\section{Repair behind the front}
\label{sec8}

As in the previous sections, we denote by $u(x,t)$ a solution of
\eqref{eq-u} whose initial data fulfill the conditions \eqref{cond-},
\eqref{cond+}, where $\delta \le \min(\delta_0, \delta_1)/2$ is small
enough so that \eqref{invasion_cond} holds. We know from
Corollary~\ref{loc-conv} that $u(x,t)$ converges to a travelling front
uniformly for $x \in (\overline x(t)-L,+\infty)$, for any $L > 0$. 
To conclude the proof of Theorem~\ref{mainthm}, it remains
to prove that $u(x,t)$ converges uniformly to $1$ far behind the
invasion point. Following again the ideas introduced in \cite{Ri}, we
shall do this using a suitable energy estimate in the {\em laboratory
frame}.

\begin{prop}\label{uniform}
There exists a sequence $t_n \to +\infty$ such that
\begin{align}\label{unifconv}
  \|u(\overline x(t_n)+\cdot,t_n) - v_*\|_{H^1_\ul} + 
  \|\dot u(\overline x(t_n)+\cdot,t_n) + s_* v_*'\|_{L^2_\ul}
  \,\xrightarrow[n\to \infty]{}\, 0~,
\end{align}
where $v_*$ is as in Corollary~\ref{loc-conv}.
\end{prop}

\begin{rem}
Using an additional argument as in \cite[Section~9.6]{Ri}, one 
can show that \eqref{unifconv} holds in fact for {\em all} sequences 
$t_n \to +\infty$. In our case, this follows from the local 
stability of the travelling front which will be established
in the last section.
\end{rem}

\begin{demo}
We recall that the solution of \eqref{eq-u} has been decomposed 
as $u(x,t) = v(x,t) + r(x,t)$, where the remainder $(r,\dot r)$ 
converges exponentially to zero as $t \to +\infty$ in the uniformly
local energy space $X = H^1_\ul(\Rm) \times L^2_\ul(\Rm)$. Using 
this remark and Corollary~\ref{loc-conv}, we can construct a 
sequence of times $\{t_n\}$ satisfying $t_{n+1} \ge t_n + n + 1$ 
for all $n \in \Nm$ and such that, for all $t \ge t_n$, 
\begin{align}\label{tndef}
  \sup_{z \ge -2n} \int_z^{z+1} \Bigl(
  |\dot u(\overline x(t)+x,t) + s_* v_*'(x)|^2 &+ 
  |u'(\overline x(t)+x,t) - v_*'(x)|^2 \\ \nonumber &+ 
  |u(\overline x(t)+x,t) - v_*(x)|^2\Bigr)\dd x 
  \,\le\, \frac{1}{n+1}~.
\end{align}
Without loss of generality, we also assume that $t_0 \ge 1$. 
Let $\theta : \Rm \to [0,1]$ be a smooth, nondecreasing function 
satisfying $\theta(x) = 0$ for $x \le -1$, $\theta(x) = 1$ for $x 
\ge 1$, and $\int_{-1}^1 \theta(x)\dd x = 1$. We define a smooth 
map $x_+ : [t_0,+\infty) \to \Rm$ in the following way. For all 
$n \in \Nm$, we set
\[
  x_+(t) \,=\, \int_0^\infty \theta'(t-\tau)\overline x(\tau)\dd \tau
  - n - \theta\Bigl(\frac{2t-t_n-t_{n+1}}{t_{n+1}-t_n}\Bigr)~,
  \quad \hbox{if } t \in [t_n,t_{n+1}]~.
\]
We recall that $t \mapsto \overline x(t)$ is upper semi-continuous, 
hence measurable. Since $\int_\Rm \theta'(x)\dd x = 1$ and 
$\int_\Rm x\theta'(x)\dd x = 0$, we have for all $t \ge 1$:
\[
  \overline x(t) - \int_0^\infty \theta'(t-\tau)\overline x(\tau)\dd \tau
  \,=\, \int_0^\infty \theta'(t-\tau)\Bigl(\overline x(t) - 
  \overline x(\tau)-s_*(t-\tau)\Bigr)\dd \tau~,
\]
and the right-hand side converges to zero as $t \to +\infty$ 
by Lemma~\ref{varxbar}. Thus, if $n \in \Nm$ is sufficiently large,
we see that 
\begin{equation}\label{xplusbdd}
  \overline x(t) -n-2 \,\le\, x_+(t) \,\le\, \overline x(t) -n+1~,
  \quad \hbox{for } t \in [t_n,t_{n+1}]~.
\end{equation}
Similarly, since $\int_\Rm \theta''(x)\dd x = 0$ and $\int_\Rm x
\theta''(x)\dd x = -1$, we have
\[
  s_* - \int_0^\infty \theta''(t-\tau)\overline x(\tau)\dd \tau
  \,=\, \int_0^\infty \theta''(t-\tau)\Bigl(\overline x(t) - 
  \overline x(\tau) -s_*(t-\tau)\Bigr)\dd \tau  
  \,\xrightarrow[t\to +\infty]{}\, 0~,
\]
hence $|x_+'(t) - s_*| \le 1$ if $t \ge 0$ is sufficiently large.
 
On the other hand, using the assumption \eqref{cond-} on the 
initial data and proceeding exactly as in the proof of 
Proposition~\ref{pinching}, we see that there exists $\xi_1 \in 
\Rm$ such that
\begin{equation}\label{xminusbdd}
  \sup_{z \le \xi_1-t/\sqrt{\alpha}} \int_z^{z+1} \Bigl(
  |\dot u(x,t)|^2 + |u'(x,t)|^2 + |u(x,t) - 1|^2\Bigr)
  \dd x \,\le\, K_1\delta_1\,e^{-\mu_1 t}~.
\end{equation}
For all $t \ge 0$, we set $x_-(t) = \xi_1-2t/\sqrt{\alpha}$.
Without loss of generality, we can assume that $x_-(t) \le x_+(t)$
for all $t \ge t_0$. 

We next define, for all $t \ge t_0$, 
\[
  \Phi(t) \,=\, \int_{\Rm} \phi(x,t) \Bigl(\frac{\alpha}{2}
  |\dot u(x,t)|^2 + \frac12|u'(x,t)|^2 + \overline V(u(x,t))\Bigr)
  \dd x~,
\]
where $\overline V(u) = V(u) - V(1) \ge 0$ and 
\[
  \phi(x,t) \,=\, \left\{\begin{array}{ccl}
  e^{x-x_-(t)} & \hbox{if} & x \le x_-(t) \\
  1 & \hbox{if} & x_-(t) \le x \le x_+(t) \\
  e^{x_+(t)-x} & \hbox{if} & x \ge x_+(t)
  \end{array}~.\right.
\]
A direct calculation shows that 
\begin{align*}
  \Phi'(t) \,=\, -\int_{\Rm} \phi(x,t) |\dot u(x,t)|^2\dd x 
  &-\int_{-\infty}^{x_-(t)} \phi \Bigl\{x_-'(t) \Bigl(
  \frac{\alpha}{2} |\dot u|^2 + \frac12|u'|^2 + \overline V(u)
  \Bigr) + \dot u u'\Bigr\}\dd x \\
  &+\int_{x_+(t)}^{\infty} \phi \Bigl\{x_+'(t) \Bigl(
  \frac{\alpha}{2} |\dot u|^2 + \frac12|u'|^2 + \overline V(u)
  \Bigr) + \dot u u'\Bigr\}\dd x~.
\end{align*}
As is clear from \eqref{xminusbdd}, the second integral in the
right-hand side converges to zero as $t \to +\infty$. Since $x_+(t) 
- \overline x(t) \to -\infty$ by \eqref{xplusbdd}, this is also
the case for the last integral. Indeed, as $v_*(x) \to 1$ when 
$x \to -\infty$, it follows from \eqref{tndef}, \eqref{xplusbdd}
that $(u,\dot u)(x_+(t)+\cdot,t)$ converges to $(1,0)$ in 
$H^1_\loc(\Rm) \times L^2_\loc(\Rm)$ as $t \to +\infty$. Since 
$\Phi(t) \ge 0$ for all times, we conclude that, given any 
$T > 0$, there exists a sequence $t_n' \to +\infty$ such that 
\begin{equation}\label{repairdissip}
  \int_{t_n'-T}^{t_n'+T} \int_\Rm \phi(x,t) |\dot u(x,t)|^2 
  \dd x \dd t \,\xrightarrow[n\to \infty]{}\, 0~.
\end{equation}

Now, we claim that 
\begin{equation}\label{repairconv}
  \sup_{z \in [x_-(t_n'),x_+(t_n')]} \int_z^{z+1} \Bigl(
  |\dot u(x,t_n')|^2 + |u'(x,t_n')|^2 + |u(x,t_n') - 1|^2\Bigr)
  \dd x \,\xrightarrow[n\to \infty]{}\, 0~.
\end{equation}
Assume on the contrary that, after extracting a subsequence, 
the left-hand side of \eqref{repairconv} is bounded from below
for all $n \in \Nm$ by some $\epsilon > 0$. Then, given any 
$\epsilon' \in (0,\epsilon)$, we can find a sequence $\{z_n\}$ 
such that $z_n \in [x_-(t_n'),x_+(t_n')]$ for all $n \in \Nm$ 
and 
\[
  \int_0^1 \Bigl(|\dot u(z_n+x,t_n')|^2 + |u'(z_n+x,t_n')|^2 + 
  |u(z_n+x,t_n') - 1|^2\Bigr)\dd x \,=\, \epsilon'~.
\]
Moreover, using \eqref{tndef}, \eqref{xminusbdd}, and the definitions
of $x_-(t)$, $x_+(t)$, we see that $x_+(t_n') - z_n \to +\infty$ and
$z_n - x_-(t_n') \to +\infty$ as $n \to \infty$. Without loss 
of generality, we can assume that $\epsilon' > 0$ is sufficiently
small so that the following property holds: if $w : \Rm \to \Rm$
is a bounded solution of the differential equation $w'' - V'(w) = 0$
such that $|w(0)-1|^2 \le 2\epsilon'$, then $w(x) = 1$ for all 
$x \in \Rm$. The existence of such an $\epsilon'$ follows 
from our assumptions \eqref{coercive}--\eqref{critval} on the 
potential $V$. 

Using once again Proposition~\ref{propcompact}, we can assume that, 
after extracting a subsequence, the sequence of functions 
$(u,\dot u)(z_n+\cdot,t_n'+\cdot)$ converges in the space 
$\CC^0([-T,T],H^1_\loc(\Rm) \times L^2_\loc(\Rm))$ to some limit 
$(w,\dot w)$ which satisfies \eqref{eq-u}. In view of 
\eqref{repairdissip}, we have $\dot w \equiv 0$, so that 
$w : \Rm \to \Rm$ is a bounded solution of the differential equation 
$w'' - V'(w) = 0$. Moreover, 
\[
  |w(0)-1|^2 \,\le\, 2\int_0^1(|w'(x)|^2 + |w(x)-1|^2)\dd x 
  \,=\, 2\epsilon'~,
\]
hence our assumption on $\epsilon'$ implies that $w \equiv 1$, 
which is clearly absurd. Thus \eqref{repairconv} is established, 
and using in addition \eqref{tndef}, \eqref{xminusbdd} we see 
that \eqref{unifconv} holds for the sequence $\{t_n'\}$. This 
concludes the proof. 
\end{demo}


\section{Local stability of the travelling front}
\label{sec9}

The aim of this final section is to show that the family \eqref{TWdef}
of travelling fronts of \eqref{eq-u} is {\em asymptotically stable with
 shift} in the uniformly local energy space $X = H^1_\ul(\Rm) \times
L^2_\ul(\Rm)$.  Together with Proposition~\ref{uniform}, this will
conclude the proof of Theorem~\ref{mainthm}. Whereas a lot is known
about local stability of travelling fronts in parabolic systems (see
e.g. \cite{Sa}), for the hyperbolic equation \eqref{eq-u} with a
bistable potential we are only aware of the note \cite{Ga} where local
stability in the usual energy space $H^1(\Rm) \times L^2(\Rm)$ is
briefly discussed.

From now on, we fix $c = c_*$ and we denote by $h : \Rm \to 
(0,1)$ the unique solution of \eqref{hdef} such that $h(0) = 
\epsilon_0$. Linearizing \eqref{eq-uc} at the steady state 
$u_c = h$, we obtain the evolution equation
\begin{equation}\label{ulin}
  \alpha u_{tt} + u_t - 2\alpha c u_{yt} \,=\, u_{yy} + c u_y - 
  g(y)u~,
\end{equation}
where $g(y) = V''(h(y))$. Proceeding as in Section~\ref{sec2}, 
it is straightforward to verify that \eqref{ulin} defines a 
$C_0$-group $\{S(t)\}_{t \in \Rm}$ of bounded linear operators 
in $X$, the generator of which is the linear operator $\AA$ 
given by 
\begin{equation}\label{AAdef}
  D(\AA) \,=\, Y~, \qquad \AA \,=\, \frac{1}\alpha 
  \begin{pmatrix} 0 && \alpha \\ \partial_y^2 + c\partial_y - g(y) 
  && -1 + 2\alpha c \partial_y \end{pmatrix}~,
\end{equation}
where $Y = H^2_\ul(\Rm) \times H^1_\ul(\Rm)$. By translation 
invariance, $\lambda = 0$ is an eigenvalue of $\AA$ with eigenfunction 
$(h',0)$. This eigenvalue is in fact simple, and the corresponding 
spectral projection reads
\begin{equation}\label{projec}
  \Pi \begin{pmatrix} u \\ v\end{pmatrix} \,=\, N 
  \begin{pmatrix} h' \\ 0\end{pmatrix}
  \int_\Rm (\psi_1 u + \psi_2 v)\dd y~,
\end{equation}
where $\psi_2(y) = e^{cy}h'(y)$, $\psi_1 = \alpha^{-1}\psi_2 
+ 2c \psi_2'$, and $N > 0$ is a normalization factor. One can 
check that $\psi_1,\psi_2$ decay exponentially to zero as 
$|y| \to \infty$. The main result of this section is:
 
\begin{prop}\label{Sspec}
There exist positive constants $C_0$ and $\nu$ such that 
\begin{equation}\label{Sestim}
  \|S(t)({\bf 1}-\Pi)\|_{\LL(X)} \,\le\, C_0\,e^{-\nu t}~, 
  \quad \hbox{for all } t \ge 0~.  
\end{equation}
\end{prop}

Using Proposition~\ref{Sspec} and classical arguments which can be
found in \cite{Sa} or \cite{He}, it is easy to show that the family of
all translates of the steady state $(h,0)$ is {\em normally
  hyperbolic} and {\em asymptotically stable} for the evolution
defined by \eqref{eq-uc} on $X$. In other words, any solution of
\eqref{eq-uc} whose initial data lie in a small tubular neighborhood
of this family of equilibria (in the topology of $X$) converges
exponentially fast as $t \to +\infty$ to some element of the family.
Now, if $u(x,t)$ is a solution of \eqref{eq-u} satisfying the
assumptions of Theorem~\ref{mainthm}, Proposition~\ref{uniform} shows
that the corresponding solution of \eqref{eq-uc} eventually enters
such a tubular neighborhood. Thus there exists $x_0 \in \Rm$ such that
\begin{equation}\label{convul}
  \|u_c(\cdot,t) - h(\cdot - x_0)\|_{H^1_\ul} + 
  \|\dot u_c(\cdot,t)\|_{L^2_\ul} \,=\, \OO(e^{-\nu t})~, 
  \quad \hbox{as } t \to +\infty~,
\end{equation}
which implies \eqref{convres} and concludes the proof of 
Theorem~\ref{mainthm}. 

It remains to prove Proposition~\ref{Sspec}. Let
\[
  g_\infty(y) \,=\, V''(1) \,\frac{1-\tanh(y)}{2} + 
  V''(0) \,\frac{1+\tanh(y)}{2}~, \quad y \in \Rm~.
\]
Obviously $g_\infty(y) \ge m$ for all $y \in \Rm$, where $m
=\min(V''(0),V''(1)) > 0$. Let $\AA_\infty$ be the linear operator on
$X$ obtained by replacing $g$ with $g_\infty$ in the definition
\eqref{AAdef} of $\AA$, and let $S_\infty(t)$ be the $C_0$-group
generated by $\AA_\infty$. For any $u \in H^1_\ul(\Rm)$, the map
$y\mapsto (g(y)-g_\infty(y))u(y)$ belongs to $H^1(\Rm)$ and converges
exponentially to zero as $y \to \pm\infty$. In particular, the linear
operator $\AA-\AA_\infty:~(u,v)\mapsto \alpha^{-1}(0, (g_\infty-g)u)$
is {\em compact} in $X$. Moreover, the group $S_\infty(t)$ is bounded
in $H^2(\Rm) \times H^1(\Rm)$ and, due to the finite speed of
propagation, it preserves the exponential decay in space.  Therefore,
the Duhamel perturbation formula
\[
  S(t)w \,=\, S_\infty(t)w + \int_0^t S_\infty(t-\tau)
  (\AA-\AA_\infty)S(\tau)w\dd \tau~, \quad w \in X~,
\]
shows that $S(t)-S_\infty(t)$ is compact in $X$ for any $t \in \Rm$.
In particular, $S(t)$ and $S_\infty(t)$ have the same essential
spectrum (in what follows, we use the notion of essential spectrum
adopted in \cite{He}). The first step in the proof of
Proposition~\ref{Sspec} is:

\begin{lemma}\label{Sinfty}
There exist positive constants $C_1$ and $\nu_\infty$ such that 
\[
  \|S_\infty(t)\|_{\LL(X)} \,\le\, C_1\,e^{-\nu_\infty t}~, 
  \quad \hbox{for all } t \ge 0~.  
\]
\end{lemma}

\begin{demo}
We proceed as in the proof of Proposition~\ref{propexist}.  Let
$\rho(y) = \exp(-\kappa |y|)$, where $\kappa > 0$ is small enough so
that $8\kappa \alpha c \le 1$, $3\kappa \le c$, and $4\kappa c \le m$.
If $u(x,t)$ is a solution of the evolution equation $\alpha u_{tt} +
u_t - 2\alpha c u_{yt} = u_{yy} + c u_y - g_\infty(y)u$ associated
with $S_\infty(t)$, we define for all $\xi \in \Rm$ and all $t \ge 0$:
\[
  \hat\EE(\xi,t) \,=\, \int_\Rm \rho(y-\xi)\Bigl(\alpha^2 
  u_t^2 + \alpha u_y^2 + \alpha g_\infty u^2 + \frac12 u^2 
  + \alpha u u_t\Bigr)(y,t)\dd y~.
\]
Since $|\alpha u u_t| \le \frac12(u^2+\alpha^2 u_t^2)$, there exists 
a constant $C_1 \ge 1$ such that
\begin{equation}\label{EEequiv}
  C_1^{-1}\|(u(\cdot,t),\dot u(\cdot,t))\|_X^2 \,\le\, 
  \sup_{\xi \in \Rm} \hat \EE(\xi,t) \,\le\, C_1 \|(u(\cdot,t),
  \dot u(\cdot,t))\|_X^2~.
\end{equation}
On the other hand, differentiating $\hat \EE(\xi,t)$ with respect
to $t$ and using our assumptions on $\kappa$, we obtain
\begin{align*}
  \partial_t \hat \EE(\xi,t) \,=\,&  -\int_\Rm \rho(y-\xi)
  (\alpha u_t^2 + u_y^2 + g_\infty u^2)(y,t)\dd y \\
  & - \int_\Rm \rho'(y-\xi)\Bigl(2\alpha u_y u_t + 2\alpha^2 c 
  u_t^2 + 2\alpha c u u_t + uu_y + \frac{c}{2}u^2\Bigr)(y,t)\dd y\\
  \,\le\,& -\frac12 \int_\Rm \rho(y-\xi) (\alpha u_t^2 + u_y^2 
  + g_\infty u^2)(y,t)\dd y \,\le\, -2\nu_\infty \hat \EE(\xi,t)~,
\end{align*}
for some $\nu_\infty > 0$. Thus $\hat \EE(\xi,t) \le \hat \EE(\xi,0)\,
e^{-2\nu_\infty t}$ for all $t \ge 0$, and using in addition 
\eqref{EEequiv} we obtain the desired estimate.
\end{demo}

Since $S(t)$ is a compact perturbation of $S_\infty(t)$ 
for any $t \ge 0$, it follows from Lemma~\ref{Sinfty} that 
the spectrum of $S(t)$ outside the disk $\{z \in \Cm \,|\,
|z| \le e^{-\nu_\infty t}\}$ consists of isolated eigenvalues with 
finite multiplicities. By the spectral mapping theorem, any
such eigenvalue has the form $e^{\lambda t}$, where $\lambda \in
\Cm$ is an eigenvalue of the generator $\AA$. If $w = (u,v) 
\in Y$ satisfies $\AA w = \lambda w$, it follows from \eqref{AAdef}
that $v = \lambda u$ and
\begin{equation}\label{nonlineig}
  u'' + c(1+2\alpha\lambda)u' - g(y)u \,=\, \lambda(1+\alpha\lambda)u~.
\end{equation}
It remains to study the nonlinear eigenvalue problem
\eqref{nonlineig}. Let
\[
  \mu_\alpha \,=\, \frac{1}{2\alpha}\Bigl(-1 + \Re 
  \sqrt{1-4\alpha m}\Bigr) \,<\, 0~,
\]
where $m = \min(V''(0),V''(1)) > 0$. One can check that 
$\mu_\alpha = \sup\{\Re(\lambda)\,|\,\lambda \in \sigma_\ess(\AA)\}$. 
The key observation is:

\begin{prop}\label{lambda}
The spectrum of $\AA$ in the half-plane
$\{z\in\Cm\,|\,\Re(z)>\mu_\alpha\}$ consists of a finite sequence of
simple real eigenvalues $0 = \lambda_0 > \lambda_1 > \dots > 
\lambda_k > \mu_\alpha$, where $k \in \Nm$ depends on $\alpha$
and $V$.
\end{prop}

The proof of Proposition~\ref{lambda} relies on the following 
elementary result: 

\begin{lemma}\label{elem}
Assume that $z \in \Cm$ satisfies $\Re(z) > \mu_\alpha$, 
and let
\[
  \gamma \,=\, \frac{c}{2}(1+2\alpha z)~, \quad 
  \delta \,=\, m + z(1+\alpha z)~.
\]
Then $\Re\sqrt{\gamma^2+\delta} > \Re(\gamma) > 0$. 
\end{lemma}

\begin{demo}
Since $\Re(z) > \mu_\alpha \ge -1/(2\alpha)$, it is clear 
that $\Re(\gamma) > 0$. On the other hand, if $z > 0$, 
we have $\gamma > 0$ and $\delta > 0$, hence $\sqrt{\gamma^2+\delta} 
> \gamma$. Thus, to prove Lemma~\ref{elem}, it is sufficient to
verify that the equality $\Re\sqrt{\gamma^2+\delta} = \Re(\gamma)$
cannot occur if $\Re(z) > \mu_\alpha$. Assume on the 
contrary that $\sqrt{\gamma^2+\delta} = \gamma +i\beta$ for some
$\beta \in \Rm$. Then $\delta = 2i\gamma\beta - \beta^2$, and 
if we set $z = z_1 + iz_2$ with $z_1, 
z_2 \in \Rm$ we obtain the relation
\begin{equation}\label{rel}
  m + z_1 + \alpha(z_1^2-z_2^2) + iz_2
  (1+2\alpha z_1) \,=\, -2\alpha\beta c z_2 -\beta^2
  + i\beta c(1+2\alpha z_1)~.
\end{equation}
Taking the imaginary part of both sides, we find $z_2 = 
\beta c$, because $1+2\alpha z_1 > 0$ by assumption. 
Using this information and taking the real part of \eqref{rel}, 
we arrive at
\begin{equation}\label{ellipse}
  \alpha\Bigl(z_1+\frac{1}{2\alpha}\Bigr)^2 + 
  \Bigl(\alpha + \frac{1}{c^2}\Bigr)^2z_2^2 \,=\, 
  \frac{1}{4\alpha} - m~.
\end{equation}
This is clearly impossible if $4\alpha m > 1$. In the converse
case, equation \eqref{ellipse} defines an ellipse in $\Cm$ 
which is entirely contained in the half-plane $\{\Re(z) 
\le \mu_\alpha\}$, thus contradicting our assumption on $z$.
\end{demo}

\medskip\noindent 
{\textbf{Proof of Proposition~\ref{lambda}:}} Assume that $\lambda
\in \sigma(\AA)$ satisfies $\Re(\lambda) > \mu_\alpha$. Then 
$\lambda$ is an eigenvalue of $\AA$, and there exists a nonzero 
$u \in H^2_\ul(\Rm)$ satisfying \eqref{nonlineig}. Since $g(y)$ 
converges exponentially to $V''(0)$ as $y \to +\infty$, we know 
that $u(y) = A \phi_1(y) + B \phi_2(y)$ for some $A,B \in \Cm$,
where $\phi_1, \phi_2$ are particular solutions of 
\eqref{nonlineig} satisfying
\[
  \lim_{y \to +\infty} \phi_1(y)\, e^{\gamma y} e^{\sqrt{\gamma^2+
  \delta_+}\,y} \,=\, 1~, \quad
  \lim_{y \to +\infty} \phi_2(y)\, e^{\gamma y} e^{-\sqrt{\gamma^2+
  \delta_+}\,y} \,=\, 1~,
\]
where $\gamma = \frac{c}{2}(1+2\alpha\lambda)$ and $\delta_+ = 
V''(0) + \lambda(1+\alpha \lambda)$, see \cite[Section 3.8]{CL}. 
But Lemma~\ref{elem} implies that $\Re\sqrt{\gamma^2+\delta_+} >
\Re(\gamma) > 0$, hence we must have $B = 0$ because $\phi_2(y)$
is unbounded as $y \to +\infty$. Thus 
\[
  u(y) \,=\, A \phi_1(y) \,\approx\, A\,e^{-\gamma y} 
  \,e^{-\sqrt{\gamma^2+\delta_+}\,y}~, \quad \hbox{as } y \to 
  +\infty~,
\]
and a similar argument shows that 
\[
  u(y) \,\approx\, C\,e^{-\gamma y} \,e^{\sqrt{\gamma^2+\delta_-}\,y}
  ~, \quad \hbox{as } y \to -\infty~,
\]
for some $C \in \Cm$, where $\delta_- = V''(1) + \lambda(1+\alpha 
\lambda)$. These obervations reveal in particular that the bounded 
solutions of \eqref{nonlineig} form a one-dimensional family, 
hence $\lambda$ is a simple eigenvalue of $\AA$.

Moreover, if we set $U(y) = e^{\gamma y}u(y)$, then $U(y)$ decays 
exponentially to zero as $y \to \pm\infty$, and a direct calculation 
shows that $U$ solves the differential equation
\begin{equation}\label{selfadjoint}
  U'' - \Bigl(g(y) + \frac{c^2}{4}\Bigr)U \,=\, 
  \lambda(1+\alpha\lambda)(1+\alpha c^2)U~, 
  \quad y \in \Rm~.
\end{equation}
Thus $\mu = \lambda(1+\alpha\lambda)(1+\alpha c^2)$ is an 
eigenvalue of the selfadjoint operator $\LL = \partial_y^2 -
(g+c^2/4)$. In particular, we have $\mu \in \Rm$, hence 
$\lambda \in \Rm$ because $\Re(\lambda) > -1/(2\alpha)$. Furthermore, 
since $\mu = 0$ is an eigenvalue of $\LL$ with eigenfunction 
$U(y) = e^{\gamma y}h'(y) < 0$, we know from Sturm-Liouville theory 
that all the other eigenvalues of $\LL$ are strictly negative. 
Finally, it follows from Bargmann's bound and the min-max principle
that $\LL$ has only a finite number of eigenvalues, see e.g. 
\cite{Si}. We conclude that the spectrum of $\AA$ in the 
half-plane $\{z\in\Cm\,|\,\Re(z)>\mu_\alpha\}$ consists of 
the eigenvalue $\lambda = 0$ and, possibly, of a finite 
number of negative eigenvalues. \hfill$\square$

\medskip It is now easy to conclude the proof of
Proposition~\ref{Sspec}.  We know from Proposition~\ref{lambda} that
$\lambda = 0$ is a simple, isolated eigenvalue of $\AA$, and that the
rest of the spectrum lies in the half-plane $\{\Re(\lambda) \le
-\hat\nu\}$ for some $\hat\nu > 0$. Going back to the semigroup
$S(t)$ generated by $\AA$, we infer that for any $t > 0$ the spectrum
of $S(t)$ is entirely contained in the disk $\{z \in \Cm\,|\, 
|z| \le e^{-\nu t}\}$, where $\nu = \min(\nu_\infty,\hat\nu)$ and
$\nu_\infty$ is as in Lemma~\ref{Sinfty}, except for the simple
eigenvalue $z = 1$ which is due to the translation invariance.
If we remove that eigenvalue by restricting $S(t)$ to the invariant
subspace $\ker \Pi$, where $\Pi$ is the spectral projection
\eqref{projec}, we obtain estimate \eqref{Sestim}, possibly 
with a slightly smaller $\nu$. The proof of Theorem~\ref{mainthm} 
is now complete. 


\end{document}